\documentclass{icmart}

\contact[matdw@bristol.ac.uk]{School of Mathematics, University of Bristol, University Walk, Clifton, Bristol BS8 1TW, United Kingdom}

\newtheorem{theorem}{Theorem}[section]

\newtheorem{conjecture}[theorem]{Conjecture}

\theoremstyle{definition}

\numberwithin{equation}{section}

\newcommand{\mmod}[1]{\,\,(\text{mod}\,\,#1)}
\def\alp{{\alpha}} \def\bfalp{{\boldsymbol \alpha}} \def\bet{{\beta}} \def\bfbet{{\boldsymbol \beta}} 
  \def\eps{{\varepsilon}}
\def\tet{{\theta}} 
\def\bfxi{{\boldsymbol \xi}}
\def\dbC{{\mathbb C}} \def\dbF{{\mathbb F}} \def\dbN{{\mathbb N}} 
\def\dbQ{{\mathbb Q}} \def\dbR{{\mathbb R}} \def\dbZ{{\mathbb Z}}
\def\grB{{\mathfrak B}} \def\grC{{\mathfrak C}}
\def\grf{{\mathfrak f}} \def\grF{{\mathfrak F}}
\def\grm{{\mathfrak m}} \def\grM{{\mathfrak M}}
\def\grS{{\mathfrak S}}
\def\Gtil{{\widetilde G}}
\def\diff{{\,{\rm d}}} 
\def\llbracket{\lbrack\;\!\!\lbrack} \def\rrbracket{\rbrack\;\!\!\rbrack}

\title[Translation invariance, exponential sums, and Waring's problem]
{Translation invariance, exponential sums, and Waring's problem}

\author[Trevor D. Wooley]
{Trevor D. Wooley}

\begin{document}

\begin{abstract} We describe mean value estimates for exponential sums of degree exceeding $2$ that 
approach those conjectured to be best possible. The vehicle for this recent progress is the 
\textit{efficient congruencing method}, which iteratively exploits the translation invariance of associated 
systems of Diophantine equations to derive powerful congruence constraints on the underlying variables. 
There are applications to Weyl sums, the distribution of polynomials modulo $1$, and other Diophantine 
problems such as Waring's problem.
\end{abstract}



\maketitle

\section{Introduction} Although pivotal to the development of vast swathes of analytic number theory in 
the twentieth century, the differencing methods devised by Weyl \cite{Wey1916} and van der Corput 
\cite{vdC1922} are in many respects unsatisfactory. In particular, they improve on the trivial estimate for 
an exponential sum by a margin exponentially small in terms of its degree. The method introduced by 
Vinogradov \cite{Vin1935, Vin1947} in 1935, based on mean values, is rightly celebrated as a great leap 
forward, replacing this exponentially weak margin by one polynomial in the degree. Nonetheless, 
Vinogradov's methods yield bounds removed from the sharpest conjectured to hold by a margin at least 
logarithmic in the degree, a defect that has endured for six decades since the era in which these ideas were 
comprehensively analysed. In this report, we describe progress since 2010 that eliminates this defect, 
placing us within a whisker of establishing in full the main conjecture of the subject.\par 

When $k,s\in \dbN$ and $\bfalp\in \dbR^k$, consider the exponential sum
\begin{equation}\label{wool1.1}
f_k(\bfalp;X)=\sum_{1\le x\le X}e(\alp_1x+\ldots +\alp_kx^k)
\end{equation}
and the mean value
\begin{equation}\label{wool1.2}
J_{s,k}(X)=\oint |f_k(\bfalp;X)|^{2s}\diff\bfalp . 
\end{equation}
Here, as usual, we write $e(z)$ for $e^{2\pi iz}$. Also, to save clutter, when $G:[0,1)^k\rightarrow \dbC$ 
is integrable, we write $\oint G(\bfalp)\diff \bfalp =\int_{[0,1]^k}G(\bfalp )\diff \bfalp $. By orthogonality, 
one sees that $J_{s,k}(X)$ counts the number of integral solutions of the system of equations
\begin{equation}\label{wool1.3}
x_1^j+\ldots +x_s^j=y_1^j+\ldots +y_s^j\quad (1\le j\le k),
\end{equation}
with $1\le x_i,y_i\le X$ $(1\le i\le s)$. Upper bounds for $J_{s,k}(X)$ are known collectively as 
{\it Vinogradov's mean value theorem}. We now focus discussion by recording the classical version of this 
theorem that emerged from the first half-century of refinements following Vinogradov's seminal paper 
\cite{Vin1935} (see in particular \cite{Hua1949, Lin1943, Vin1947}), culminating in the papers of 
Karatsuba \cite{Kar1973} and Stechkin \cite{Ste1975}.

\begin{theorem}\label{wooltheorem1.1} There is an absolute constant $A>0$ having the property that, 
whenever $s$, $r$ and $k$ are natural numbers with $s\ge rk$, then
\begin{equation}\label{wool1.4}
J_{s,k}(X)\le C(k,r)X^{2s-k(k+1)/2+\Delta_{s,k}},
\end{equation}
where $\Delta_{s,k}=\tfrac{1}{2}k^2(1-1/k)^r$ and $C(k,r)=\min\{ k^{Ask},k^{Ak^3}\}$.
\end{theorem}

We will not concern ourselves with the dependence on $s$ and $k$ of constants such as $C(k,r)$ 
appearing in bounds for $J_{s,k}(X)$ and its allies (but see \cite{Woo1993} for improvements in this 
direction). Although significant in applications to the zero-free region of the Riemann zeta function, this is 
not relevant to those central to this paper. Thus, implicit constants in the notation of Landau and 
Vinogradov will depend at most on $s$, $k$ and $\eps$, unless otherwise indicated\footnote{Given a 
complex-valued function $f(t)$ and positive function $g(t)$, we use Vinogradov's notation $f(t)\ll g(t)$, or 
Landau's notation $f(t)=O(g(t))$, to mean that there is a positive number $C$ for which $f(t)\le Cg(t)$ for 
all large enough values of $t$. Also, we write $f(t)\gg g(t)$ when $g(t)\ll f(t)$. If $C$ depends on certain 
parameters, then we indicate this by appending these as subscripts to the notation. Also, we write 
$f(t)=o(g(t))$ when $f(t)/g(t)\rightarrow 0$ as $t\rightarrow \infty$. Finally, we use the convention that 
whenever $\eps$ occurs in a statement, then the statement holds for each fixed $\eps>0$.}.\par

When $k\ge 2$, the exponent $\Delta_{s,k}$ of Theorem \ref{wooltheorem1.1} satisfies 
$\Delta_{s,k}\le k^2e^{-s/k^2}$, and so $\Delta_{s,k}=O(1/\log k)$ for 
$s\ge k^2(2\log k+\log \log k)$. One can refine (\ref{wool1.4}) to obtain an asymptotic formula when 
$s$ is slightly larger (see \cite[Theorem 3.9]{ACK2004}, for example).

\begin{theorem}\label{wooltheorem1.2} Let $k,s\in \dbN$ and suppose that 
$s\ge k^2(2\log k+\log \log k+5)$. Then there exists a positive number $\grC(s,k)$ with 
$J_{s,k}(X)\sim \grC(s,k) X^{2s-k(k+1)/2}$.
\end{theorem}

With these theorems in hand, we consider the motivation for investigating the sums $f_k(\bfalp;X)$. 
Many number-theoretic functions may be estimated in terms of such sums. Thus, when $\text{Re}(s)$ is 
close to $1$, estimates for the Riemann zeta function $\zeta(s)$ stem from partial summation and 
Taylor expansions for $\log (1+x/N)$, since
\[ \sum_{N<n\le N+X}n^{-it}=N^{-it}\sum_{1\le x\le X}e\Bigl( -\frac{t}{2\pi}\log (1+x/N)\Bigr) .\]
On the other hand, specialisations of $f_k(\bfalp;X)$ arise naturally in applications of interest. 
Indeed, work on the asymptotic formula in Waring's problem depends on the sum obtained by 
setting $\alpha_1=\ldots =\alpha_{k-1}=0$ and $\alpha_k=\beta$, namely
\begin{equation}\label{wool1.6}
g_k(\beta;X)=\sum_{1\le x\le X}e(\beta x^k).
\end{equation}
Writing $R_{s,k}(n)$ for the number of representations of $n$ as the sum of $s$ positive integral 
$k$th powers, one finds by orthogonality that
\[ R_{s,k}(n)=\int_0^1g_k(\beta ;n^{1/k})^se(-\beta n)\diff \beta .\]

\par The uninitiated reader will wonder why one should focus on estimates for the mean value 
$J_{s,k}(X)$ when many applications depend on pointwise estimates for $f_k(\bfalp;X)$. Vinogradov 
observed that mean value estimates suffice to obtain useful pointwise estimates for $f_k(\bfalp;X)$. 
To see why this is the case, note first that $|f_k(\bfbet;X)|$ differs little from $|f_k(\bfalp;X)|$ 
provided that the latter is large, and in addition $|\bet_j-\alp_j|$ is rather smaller than $X^{-j}$ 
for each $j$, so that $\bfbet$ lies in a small neighbourhood of $\bfalp$ having measure of order 
$X^{-k(k+1)/2}$. Second, one sees from (\ref{wool1.1}) that for each integer $h$ the sum 
$f_k(\bfalp;X)$ may be rewritten in the form
\[ f_k(\bfalp;X)=\sum_{1-h\le x\le X-h}e\bigl(\alp_1(x+h)+\ldots +\alp_k(x+h)^k\bigr). \]
By estimating the tails of this sum and applying the binomial theorem to identify the coefficient of each 
monomial $x^j$, one obtains new $k$-tuples $\bfalp^{(h)}$ for which 
$f_k(\bfalp;X)=f_k(\bfalp^{(h)};X)+O(|h|)$. These ideas combine to show that one large value 
$|f_k(\bfalp;X)|$ generates a collection of neighbourhoods $\grB(h)$, with the property that whenever 
$\bfbet\in \grB(h)$, then $|f_k(\bfbet;X)|$ is almost as large as $|f_k(\bfalp;X)|$. Given $N$ disjoint 
such neighbourhoods over which to integrate $|f_k(\bfbet;X)|^{2s}$, non-trivial estimates for 
$|f_k(\bfalp;X)|$ follow from the relation $NX^{-k(k+1)/2}|f_k(\bfalp;X)|^{2s}\ll J_{s,k}(X)$. This circle 
of ideas leads to the following theorem (see \cite{Bom1990} and \cite[Theorem 5.2]{Vau1997}).

\begin{theorem}\label{wooltheorem1.3} Let $k$ be an integer with $k\ge 2$, and let $\bfalp\in \dbR^k$. 
Suppose that there exists a natural number $j$ with $2\le j\le k$ such that, for some $a\in \dbZ$ and 
$q\in \dbN$ with $(a,q)=1$, one has $|\alp_j-a/q|\le q^{-2}$. Then one has
$$f_k(\bfalp;X)\ll \left( X^{k(k-1)/2}J_{s,k-1}(2X)(q^{-1}+X^{-1}+qX^{-j})\right)^{1/(2s)}\log (2X).$$
\end{theorem} 

To illustrate the power of this theorem, suppose that $k$ is large, and $\bet$ satisfies the condition that, 
whenever $b\in \dbZ$ and $q\in \dbN$ satisfy $(b,q)=1$ and $|q\bet-b|\le X^{1-k}$, then $q>X$. By 
substituting the conclusion of Theorem \ref{wooltheorem1.1} into Theorem \ref{wooltheorem1.3}, one 
finds that $g_k(\bet;X)\ll X^{1-\sigma(k)}$, where $\sigma(k)^{-1}=(4+o(1))k^2\log k$.

\section{Translation invariance and a congruencing idea} A key feature of the system of equations 
(\ref{wool1.3}) is \emph{translation-dilation invariance}. Thus, the pair 
${\boldsymbol x},{\boldsymbol y}$ is an integral solution of the system
\begin{equation}\label{woola.1}
x_1^j+\ldots +x_t^j=y_1^j+\ldots +y_t^j\quad (1\le j\le k),
\end{equation}
if and only if, for any $\xi\in \dbZ$ and $q\in \dbN$, the pair ${\boldsymbol x},{\boldsymbol y}$ 
satisfies the system
\begin{equation}\label{woola.2}
(qx_1+\xi)^j+\ldots +(qx_t+\xi)^j=(qy_1+\xi)^j+\ldots +(qy_t+\xi)^j\quad (1\le j\le k).
\end{equation}
This property ensures that $J_t(X)$ is homogeneous with respect to restriction to arithmetic 
progressions\footnote{In this section we consider $k$ to be fixed, and hence we drop mention of $k$ 
from our notations.}. Let $M=X^\theta$ be a parameter to be chosen later, consider a set ${\cal P}$ of 
$\lceil k^2/\theta \rceil$ primes $p$ with $M<p\le 2M$, and fix some $p\in {\cal P}$. Also, define
\[ \grf_c(\bfalp;\xi)=\sum_{\substack{1\le x\le X\\ x\equiv \xi \mmod{p^c}}}
e(\alp_1x+\ldots +\alp_kx^k).\]
Since $\oint |\grf_c(\bfalp;\xi)|^{2t}\diff \bfalp$ counts the number of solutions of (\ref{woola.2}), 
with $q=p^c$, for which\footnote{Here we make use of slightly unconventional vector notation. Thus, 
we write ${\boldsymbol z}\equiv \xi\mmod{q}$ when $z_i\equiv \xi\mmod{q}$ for $1\le i\le t$, or 
$a\le {\boldsymbol z}\le b$ when $a\le z_i\le b$ $(1\le i\le t)$, and so on.} 
$(1-\xi)/q\le {\boldsymbol x},{\boldsymbol y}\le (X-\xi)/q$, by translation-dilation invariance, it 
counts solutions of (\ref{woola.1}) under the same conditions on ${\boldsymbol x}$ and 
${\boldsymbol y}$. Thus
\begin{equation}\label{woola.3}
\max_{1\le \xi\le p^c}\oint |\grf_c(\bfalp;\xi)|^{2t}\diff \bfalp \ll 1+J_t(X/M^c).
\end{equation}

Translation invariance also generates useful auxiliary congruences. Let $t=s+k$, and 
consider the solutions of (\ref{woola.1}) with $1\le {\boldsymbol x},{\boldsymbol y}\le X$. The 
number of solutions $T_0$ in which $x_i=x_j$ for some $1\le i<j\le k$ may be bounded via orthogonality 
and H\"older's inequality, giving $T_0\ll J_t(X)^{1-1/(2t)}$. Given a {\it conditioned} solution with 
$x_i\ne x_j$ for $1\le i<j\le k$, there exists a prime $p\in {\cal P}$ with $x_i\not\equiv x_j\mmod{p}$ 
for $1\le i<j\le k$. Let $\Xi_c(\xi)$ denote the set of $k$-tuples $(\xi_1,\ldots ,\xi_k)$, 
with $1\le \bfxi\le p^{c+1}$ and $\bfxi\equiv \xi\mmod{p^c}$, and satisfying the property that 
$\xi_i\not \equiv \xi_j\mmod{p^{c+1}}$ for $i\ne j$. Also, put
\[ \grF_c(\bfalp;\xi)=\sum_{\bfxi\in \Xi_c(\xi)}\grf_{c+1}(\bfalp;\xi_1)\ldots \grf_{c+1}(\bfalp;\xi_k),\]
and define
\begin{equation}\label{woola.3a}
I_{a,b}(X)=\max_{\xi,\eta}\oint |\grF_a(\bfalp;\xi)^2\grf_b(\bfalp;\eta)^{2s}|\diff \bfalp .
\end{equation}
Then for some $p\in {\cal P}$, which we now fix, the number $T_1$ of conditioned solutions satisfies
\begin{equation}\label{woola.3b}
T_1\ll \oint \grF_0(\bfalp;0)f(\bfalp;X)^sf(-\bfalp;X)^{s+k}\diff \bfalp .
\end{equation}
Thus, by Schwarz's inequality and orthogonality, one has $T_1\ll I_{0,0}(X)^{1/2}J_t(X)^{1/2}$. By 
combining the above estimates for $T_0$ and $T_1$, we derive the upper bound 
$J_t(X)\ll J_t(X)^{1-1/(2t)}+I_{0,0}(X)^{1/2}J_t(X)^{1/2}$, whence $J_t(X)\ll I_{0,0}(X)$.\par

By H\"older's inequality, one finds that
\[ |f(\bfalp;X)|^{2s}=\biggl| \sum_{\eta=1}^p \sum_{\substack{1\le x\le X\\ x\equiv \eta\mmod{p}}}
e(\alp_1x+\ldots +\alp_kx^k)\biggr|^{2s}\le p^{2s-1}\sum_{\eta=1}^p|\grf_1(\bfalp;\eta)|^{2s}.\]
Thus, on noting the trivial relation $f(\bfalp;X)=\grf_0(\bfalp;\eta)$, one sees from (\ref{woola.3a}) 
that
\begin{equation}\label{woola.a}
J_t(X)\ll I_{0,0}(X)\ll M^{2s}\max_{\xi,\eta}\oint|\grF_0(\bfalp;\xi)^2\grf_1(\bfalp;\eta)^{2s}|\diff 
\bfalp =M^{2s}I_{0,1}(X).
\end{equation}
The mean value underlying $I_{0,1}(X)$ counts the number of integral solutions of
\[ \sum_{i=1}^k(x_i^j-y_i^j)=\sum_{l=1}^s\left( (pu_l+\eta)^j-(pv_l+\eta)^j\right)\quad (1\le j\le k),
\]
with $1\le {\boldsymbol x},{\boldsymbol y}\le X$ and 
$1\le p{\boldsymbol u}+\eta,p{\boldsymbol v}+\eta\le X$, in which $x_i\not\equiv x_j\mmod{p}$ for 
$i\ne j$, and similarly for ${\boldsymbol y}$. Translation invariance leads from these equations to
\[ \sum_{i=1}^k\left( (x_i-\eta)^j-(y_i-\eta)^j\right)=p^j\sum_{l=1}^s(u_l^j-v_l^j)\quad 
(1\le j\le k),\]
and hence
\begin{equation}\label{woola.5}
(x_1-\eta)^j+\ldots +(x_k-\eta)^j\equiv (y_1-\eta)^j+\ldots +(y_k-\eta)^j\mmod{p^j}\quad 
(1\le j\le k).
\end{equation}
Since the $x_i$ are distinct modulo $p$, Hensel's lemma shows that, for each fixed choice of 
${\boldsymbol y}$, there are at most $k!p^{k(k-1)/2}$ choices for ${\boldsymbol x}\mmod{p^k}$ 
satisfying (\ref{woola.5}). An application of Cauchy's inequality shows from here that
\begin{equation}\label{woola.6}
I_{0,1}(X)\ll M^{k(k-1)/2}\max_\eta\oint \biggl( \sum_{\nu=1}^{p^k}|\grf_k(\bfalp;\nu)|^2\biggr)^k
|\grf_1(\bfalp;\eta)|^{2s}\diff \bfalp .
\end{equation}

\par Although our notation has been crafted for later discussion of efficient congruencing, the classical 
approach remains visible. One applies (\ref{woola.6}) with $\theta=1/k$, so that $p^k>X$. Thus 
$|\grf_k(\bfalp;\nu)|\le 1$, and it follows from (\ref{woola.a}) and (\ref{woola.6}) that
\[ J_t(X)\ll M^{2s}I_{0,1}(X)\ll M^{2s+k(k-1)/2}(M^k)^k
\max_\eta \oint |\grf_1(\bfalp;\eta)|^{2s}\diff \bfalp .\]
It therefore follows from (\ref{woola.3}) that
\[J_{s+k}(X)\ll M^{2s+k(k-1)/2}X^kJ_s(X/M)\ll X^{2k}(X^{1/k})^{2s-k(k+1)/2}J_s(X^{1-1/k}).\]
This iterative relation leads from the bound $J_{k,k}(X)\ll X^k$ to the estimate presented in Theorem 
\ref{wooltheorem1.1}. Early authors, such as Vinogradov and Hua, made use of short real intervals in 
place of congruences, the modern shift to congruences merely adjusting the point of view from the infinite 
place to a finite place.

\section{Lower bounds and the Main Conjecture} Write $T_s(X)$ for the number of diagonal solutions of 
(\ref{wool1.3}) with $1\le {\boldsymbol x},{\boldsymbol y}\le X$ and 
$\{x_1,\ldots ,x_s\}=\{y_1,\ldots ,y_s\}$. Then $J_{s,k}(X)\ge T_s(X)=s!X^s+O_s(X^{s-1})$. 
Meanwhile, when $1\le {\boldsymbol x},{\boldsymbol y}\le X$, one 
has $|(x_1^j-y_1^j)+\ldots +(x_s^j-y_s^j)| \le sX^j$. Hence
\[ [X]^{2s}=\sum_{|h_1|\le sX}\ldots \sum_{|h_k|\le sX^k}\oint |f_k(\bfalp;X)|^{2s}
e(-\alp_1h_1-\ldots -\alp_kh_k)\diff \bfalp ,\]
and we deduce from the triangle inequality in combination with (\ref{wool1.2}) that
\[X^{2s}\ll \sum_{|h_1|\le sX}\ldots \sum_{|h_k|\le sX^k}J_{s,k}(X)\ll X^{k(k+1)/2}J_{s,k}(X).\]
Thus we conclude that $J_{s,k}(X)\gg X^s+X^{2s-k(k+1)/2}$, a lower bound that guides a heuristic 
application of the circle method towards the following conjecture.

\begin{conjecture}[The Main Conjecture]\label{woolconj2.1}
Suppose that $s$ and $k$ are natural numbers. Then for each $\eps>0$, one has 
$J_{s,k}(X)\ll X^\eps (X^s+X^{2s-k(k+1)/2})$.
\end{conjecture}

We emphasise that the implicit constant here may depend on $\eps$, $s$ and $k$. The critical case of 
the Main Conjecture with $s=k(k+1)/2$ has special significance.

\begin{conjecture}\label{woolconj2.2}
When $k\in \dbN$ and $\eps>0$, one has $J_{k(k+1)/2,k}(X)\ll X^{k(k+1)/2+\eps}$.
\end{conjecture}

Suppose temporarily that this critical case of the Main Conjecture holds. Then, when $s\ge k(k+1)/2$, 
one may apply a trivial estimate for $f_k(\bfalp;X)$ to show that
\[ J_{s,k}(X)\le X^{2s-k(k+1)}\oint |f_k(\bfalp;X)|^{k(k+1)}\diff \bfalp \ll X^{2s-k(k+1)/2+\eps},\]
and when $s<k(k+1)/2$, one may instead apply H\"older's inequality to obtain
\[ J_{s,k}(X)\le \biggl( \oint |f_k(\bfalp;X)|^{k(k+1)}\diff \bfalp \biggl)^{\tfrac{2s}{k(k+1)}}
\ll X^{s+\eps}.\]
In both cases, therefore, the Main Conjecture is recovered from the critical case.\par

Until 2014, the critical case of the Main Conjecture was known to hold in only two cases. The case 
$k=1$ is trivial. The case $k=2$, on the other hand, depends on bounds for the number of integral 
solutions of the simultaneous equations
\begin{equation}\label{wool2.3}
\left.\begin{aligned}
x_1^2+x_2^2+x_3^2&=y_1^2+y_2^2+y_3^2\\
x_1+x_2+x_3&=y_1+y_2+y_3
\end{aligned}\right\} ,
\end{equation}
with $1\le x_i,y_i\le X$. From the identity $(a+b-c)^2-(a^2+b^2-c^2)=2(a-c)(b-c)$, one finds that the 
solutions of (\ref{wool2.3}) satisfy $(x_1-y_3)(x_2-y_3)=(y_1-x_3)(y_2-x_3)$. From here, elementary 
estimates for the divisor function convey us to the bound $J_{3,2}(X)\ll X^{3+\eps}$, so that Conjecture 
\ref{woolconj2.2} and the Main Conjecture hold when $k=2$. In fact, improving on earlier work of 
Rogovskaya \cite{Rog1986}, it was shown by Blomer and Br\"udern \cite{BB2010} that
\[ J_{3,2}(X)=\frac{18}{\pi^2}X^3\log X+\frac{3}{\pi^2}
\biggl( 12\gamma -6\frac{\zeta'(2)}{\zeta(2)}-5\biggr) X^3+O(X^{5/2}\log X).\]
In particular, the factor $X^\eps$ cannot be removed from the statements of Conjectures 
\ref{woolconj2.1} and \ref{woolconj2.2}. However, a careful heuristic analysis of the circle method 
reveals that when $(s,k)\ne (3,2)$, the Main Conjecture should hold with $\eps=0$. See 
\cite[equation (7.5)]{Vau1997} for a discussion that records precisely such a conjecture.\par

The classical picture of the Main Conjecture splits naturally into two parts: small $s$ and large $s$. 
When $1\le s\le k$, the relation $J_{s,k}(X)=T_s(X)\sim s!X^s$ is immediate from Newton's formulae 
concerning roots of polynomials. Identities analogous to that above yield multiplicative relations 
amongst variables in the system (\ref{wool1.3}) when $s=k+1$. In this way, Hua \cite{Hua1947} 
confirmed the Main Conjecture for $s\le k+1$ by obtaining the bound $J_{k+1,k}(X)\ll X^{k+1+\eps}$. 
Vaughan and Wooley have since obtained the asymptotic formula 
$J_{k+1,k}(X)=T_{k+1}(X)+O(X^{\theta_k+\eps})$, where $\tet_3=\frac{10}{3}$ 
\cite[Theorem 1.5]{VW1995} and $\tet_k=\sqrt{4k+5}$ $(k\ge 4)$ 
\cite[Theorem 1]{VW1997}. Approximations to the Main Conjecture of the type 
$J_{s,k}(X)\ll X^{s+\delta_{s,k}}$, with $\delta_{s,k}$ small, can be obtained for larger values of 
$s$. Thus, on writing $\gamma=s/k$, the work of Arkhipov and Karatsuba \cite{AK1978} shows that 
permissible exponents $\delta_{s,k}$ exist with $\delta_{s,k}\ll \gamma^{3/2}k^{1/2}$, Tyrina 
\cite{Tyr1987} gets $\delta_{s,k}\ll \gamma^2$, and Wooley \cite[Theorem 1]{Woo1994} obtains 
$\delta_{s,k}=\exp(-Ak/\gamma^2)$, when $s\le k^{3/2}(\log k)^{-1}$, for a certain positive constant 
$A$.\par

We turn next to large values of $s$. When $k\in \dbN$, denote by $H(k)$ the least integer for which 
the Main Conjecture for $J_{s,k}(X)$ holds whenever $s\ge H(k)$. Theorem \ref{wooltheorem1.2} gives 
$H(k)\le (2+o(1))k^2\log k$, a consequence of the classical estimate (\ref{wool1.4}) with permissible 
exponent $\Delta_{s,k}=k^2e^{-s/k^2}$. In 1992, the author \cite{Woo1992b} found a means of 
combining Vinogradov's methods with the {\it efficient differencing method} (see \cite{Woo1992a}, and 
the author's previous ICM lecture \cite{Woo2002a}), obtaining $\Delta_{s,k}\approx k^2e^{-2s/k^2}$. 
This yields $H(k)\le (1+o(1))k^2\log k$ (see \cite{Woo1996}), halving the previous bound. Meanwhile, 
Hua \cite[Theorem 7]{Hua1947} has applied Weyl differencing to bound $H(k)$ for small $k$. We 
summarise the classical status of the Main Conjecture in the following theorem.

\begin{theorem}\label{wootheorem2.3} The Main Conjecture holds for $J_{s,k}(X)$ when:
\begin{enumerate}
\item[(i)] $k=1$ and $2$;
\item[(ii)] $k\ge 2$ and $1\le s\le k+1$;
\item[(iii)] $s\ge H(k)$, where $H(3)=8$, $H(4)=23$, $H(5)=55$, $H(6)=120$, ..., 
and $H(k)=k^2(\log k+2\log \log k+O(1))$.
\end{enumerate}
\end{theorem}

\section{The advent of efficient congruencing} The introduction of the 
{\it efficient congruencing method} \cite{Woo2012a} at the end of 2010 has transformed our 
understanding of the Main Conjecture. Incorporating subsequent developments \cite{FW2014, Woo2013}, 
and the multigrade enhancement of the method \cite{Woo2014a, Woo2014b, Woo2014c}, we can 
summarise the current state of affairs in the form of a theorem.

\begin{theorem}\label{wooltheorem3.1} The Main Conjecture holds for $J_{s,k}(X)$ when:
\begin{enumerate}
\item[(i)] $k=1$, $2$ and $3$;
\item[(ii)] $1\le s\le D(k)$, where $D(4)=8$, $D(5)=10$, $D(6)=17$, $D(7)=20$, ..., and 
$D(k)=\frac{1}{2}k(k+1)-\frac{1}{3}k+O(k^{2/3})$;
\item[(iii)] $k\ge 3$ and $s\ge H(k)$, where $H(k)=k(k-1)$.
\end{enumerate}
\end{theorem}

As compared to the classical situation, there are three principal advances:\vskip.1cm

\noindent (a) First, the Main Conjecture holds for $J_{s,k}(X)$ in the cubic case $k=3$ (see 
\cite[Theorem 1.1]{Woo2014c}), so that $J_{s,3}(X)\ll X^\eps (X^s+X^{2s-6})$. This is the first 
occasion, for any polynomial Weyl sum of degree exceeding $2$, that the conjectural mean value estimates 
have been established in full, even if the underlying variables are restricted to lie in such special sets as the 
smooth numbers.\vskip.1cm

\noindent (b) Second, the Main Conjecture holds in the form $J_{s,k}(X)\ll X^{s+\eps}$ provided that 
$1\le s\le \tfrac{1}{2}k(k+1)-\tfrac{1}{3}k+O(k^{2/3})$, which as $k\rightarrow \infty$ represents 
$100\%$ of the critical interval $1\le s\le k(k+1)/2$ (see \cite[Theorem 1.3]{Woo2014b}). The classical 
result reported in Theorem \ref{wootheorem2.3}(ii) only provides such a conclusion for $1\le s\le k+1$, 
amounting to $0\%$ of the critical interval. Here, the first substantial advance was achieved by Ford and 
Wooley \cite[Theorem 1.1]{FW2014}, giving the Main Conjecture for 
$1\le s\le \tfrac{1}{4}(k+1)^2$. Although Theorem \ref{wooltheorem3.1}(ii) comes within 
$\left(\tfrac{1}{3}+o(1)\right)k$ variables of proving the critical case of the Main Conjecture, it seems that 
a new idea is required to replace this defect by $(c+o(1))k$, for some real number $c$ with 
$c<\tfrac{1}{3}$.\vskip.1cm 

\noindent (c) Third, the Main Conjecture holds in the form $J_{s,k}(X)\ll X^{2s-k(k+1)/2+\eps}$ 
for $s\ge k(k-1)$. The classical result reported in Theorem \ref{wootheorem2.3}(iii) provides such a 
conclusion for $s\ge (1+o(1))k^2\log k$, a constraint weaker by a factor $\log k$. So far as applications 
are concerned, this is by far the most significant advance thus far captured by the efficient congruencing 
method. The initial progress \cite[Theorem 1.1]{Woo2012a} shows that the Main Conjecture holds for 
$s\ge k(k+1)$, already within a factor $2$ of the critical exponent $s=k(k+1)/2$. Subsequently, this 
constraint was improved first to $s\ge k^2-1$, and then to $s\ge k^2-k+1$ (see 
\cite[Theorem 1.1]{Woo2013} and \cite[Corollary 1.2]{Woo2014a}). The further modest progress 
reported in Theorem \ref{wooltheorem3.1}(iii) was announced in \cite[Theorem 1.2]{Woo2014c}, and will 
appear in a forthcoming paper.\vskip.1cm 

Prior to the advent of efficient congruencing, much effort had been spent on refining estimates of the 
shape $J_{s,k}(X)\ll X^{2s-k(k+1)/2+\Delta_{s,k}}$, with the permissible exponent $\Delta_{s,k}$ as 
small as possible (see \cite{BW2012, For2002, Woo1992b, Woo1994}). Of great significance for 
applications, efficient congruencing permits substantially sharper bounds to be obtained for such exponents 
than were hitherto available. Such ideas feature in \cite[Theorem 1.4]{Woo2013}, and the discussion 
following \cite[Theorem 1.2]{FW2014} shows that when $\tfrac{1}{4}\le \alp\le 1$ and $s=\alp k^2$, 
then the exponent $\Delta_{s,k}=(1-\sqrt{\alp})^2k^2+O(k)$ is permissible. Thus, in particular, the 
critical exponent $\Delta_{k(k+1)/2,k}=(\tfrac{3}{2}-\sqrt{2})k^2$ is permissible. By combining 
\cite[Theorem 1.5]{Woo2014b} and the discussion following \cite[Corollary 1.2]{Woo2014a}, one arrives 
at the following improvement.

\begin{theorem}\label{wooltheorem3.2} When $k$ is large, there is a positive number 
$C(s)\le \tfrac{1}{3}$ for which
$$J_{s,k}(X)\ll X^{(C(s)+o(1))k}\bigl( X^s+X^{2s-k(k+1)/2}\bigr),$$
When $\alp\in [\frac{5}{8},1]$, moreover, one may take 
$C(\alp k^2)\le \bigl(2-3\alp+(2\alp-1)^{3/2}\bigr)/(3\alp)$.
\end{theorem}

We finish this section by noting that Theorem \ref{wooltheorem3.1}(iii) permits a substantial improvement 
in the conclusion of Theorem \ref{wooltheorem1.2}. 

\begin{theorem}\label{wooltheorem3.3} Let $k,s\in \dbN$ and suppose that $s\ge k^2-k+1$. Then there 
exists a positive number $\grC(s,k)$ with $J_{s,k}(X)\sim \grC(s,k) X^{2s-k(k+1)/2}$.
\end{theorem}

\section{A sketch of the efficient congruencing method} Although complicated in detail, the ideas 
underlying efficient congruencing are accessible given some simplifying assumptions. In this section, we 
consider $k$ to be fixed, and drop mention of $k$ from our notation. Let $t=(u+1)k$, 
where $u\ge k$ is an integer, and put $s=uk$. We define
\[ \lambda_t=\underset{X\rightarrow \infty}{\lim \sup}\left(\log J_t(X)\right)/(\log X).\]
Thus, for each $\eps>0$, one has the bound $J_t(X)\ll X^{\lambda_t+\eps}$. Our goal is to establish that 
$\lambda_t=2t-k(k+1)/2$, as predicted by the Main Conjecture. Define $\Lambda$ via the relation 
$\lambda_t=2t-\tfrac{1}{2}k(k+1)+\Lambda$. We suppose that $\Lambda>0$, and seek a contradiction 
in order to show that $\Lambda=0$. Our method rests on an $N$-fold iteration related to the approach of 
\S2, where $N$ is sufficiently large in terms of $u$, $k$ and $\Lambda$. Let $\theta=(16k)^{-2N}$, put 
$M=X^\theta$, and consider a prime number $p$ with $M<p\le 2M$. Also, let $\delta>0$ be small in 
terms of all these parameters, so that $8\delta<N(k/u)^N\Lambda \theta $.\par

Define the mean value
\[ K_{a,b}(X)=\max_{\xi,\eta}\oint |\grF_a(\bfalp;\xi)^2\grF_b(\bfalp;\eta)^{2u}|\diff \bfalp ,\]
and introduce the normalised mean values
\[ \llbracket K_{a,b}(X)\rrbracket =\frac{K_{a,b}(X)}{(X/M^a)^{2k-k(k+1)/2}(X/M^b)^{2s}}\quad 
\text{and}\quad \llbracket J_t(X)\rrbracket =\frac{J_t(X)}{X^{2t-k(k+1)/2}}.\]
Then whenever $X$ is sufficiently large in terms of the ambient parameters, one has 
$\llbracket J_t(X)\rrbracket >X^{\Lambda -\delta}$ and, when $X^{1/2}\le Y\le X$, we have the bound 
$\llbracket J_t(Y)\rrbracket \le Y^{\Lambda+\delta}$.\par

We begin by observing that an elaboration of the argument delivering (\ref{woola.3b}) can be fashioned 
to replace (\ref{woola.a}) with the well-conditioned relation
\[ J_t(X)\ll M^{2s}\max_{\xi,\eta}\oint |\grF_0(\bfalp;\xi)^2\grF_1(\bfalp;\eta)^{2u}|\diff 
\bfalp =M^{2s}K_{0,1}(X).\]
Here we have exercised considerable expedience in ignoring controllable error terms. Moreover, one may 
need to replace $K_{0,1}(X)$ by the surrogate $K_{0,1+h}(X)$, for a suitable integer $h$. An analogue of 
the argument leading to (\ref{woola.6}) yields the bound
\[ K_{0,1}(X)\ll M^{k(k-1)/2}\max_\eta \oint \biggl( \sum_{\nu=1}^{p^k}|\grf_k(\bfalp;\nu)|^2
\biggr)^k|\grF_1(\bfalp;\eta)|^{2u}\diff \bfalp .\]
By H\"older's inequality, one finds first that
\[ \biggl( \sum_{\nu=1}^{p^k}|\grf_k(\bfalp;\nu)|^2\biggr)^k\le (p^k)^{k-1}
\sum_{\nu=1}^{p^k}|\grf_k(\bfalp;\nu)|^{2k},\]
and then
\[ K_{0,1}(X)\ll M^{k(k-1)/2}(M^k)^k\max_{\eta,\nu}\left( 
T_1(\eta)^{1-1/u}T_2(\eta,\nu)^{1/u}\right) ,\]
where
\[ T_1(\eta)=\oint |\grF_1(\bfalp;\eta)|^{2u+2}\diff \bfalp \quad \text{and}\quad T_2(\eta,\nu)
=\oint |\grF_1(\bfalp;\eta)^2\grf_k(\bfalp;\nu)^{2s}|\diff \bfalp .\]
On considering the underlying Diophantine systems, one finds that $T_1(\eta)$ may 
be bounded via (\ref{woola.3}), while $T_2(\eta,\nu)$ may be bounded in terms of $K_{1,k}(X)$. Thus
\[ J_t(X)\ll M^{2s+k(k-1)/2}(M^k)^kJ_t(X/M)^{1-1/u}K_{1,k}(X)^{1/u}.\]
A modicum of computation therefore confirms that
\begin{equation}\label{woolb.5}
\llbracket J_t(X)\rrbracket \ll \llbracket J_t(X/M)\rrbracket^{1-1/u}\llbracket K_{1,k}(X)\rrbracket^{1/u}.
\end{equation}

\par The mean value underlying $K_{1,k}(X)$ counts the number of integral solutions of
\[ \sum_{i=1}^k(x_i^j-y_i^j)=\sum_{l=1}^s\left( (p^ku_l+\eta)^j-(p^kv_l+\eta)^j\right) \quad 
(1\le j\le k),\]
with $1\le {\boldsymbol x},{\boldsymbol y}\le X$ and 
$1\le p^k{\boldsymbol u}+\eta,p^k{\boldsymbol v}+\eta\le X$ having suitably conditioned coordinates. 
In particular, one has ${\boldsymbol x}\equiv {\boldsymbol y}\equiv \xi\mmod{p}$ but 
$x_i\not \equiv x_j\mmod{p^2}$ for $i\ne j$, and similarly for ${\boldsymbol y}$. Translation invariance 
leads from these equations to
\[ \sum_{i=1}^k\left( (x_i-\eta)^j-(y_i-\eta)^j\right) =p^{jk}\sum_{l=1}^s(u_l^j-v_l^j)\quad 
(1\le j\le k),\]
and hence to the congreunces
\begin{equation}\label{woolb.6}
(x_1-\eta)^j+\ldots +(x_k-\eta)^j\equiv (y_1-\eta)^j+\ldots +(y_k-\eta)^j\mmod{p^{jk}}
\quad (1\le j\le k).
\end{equation}
Since the $x_i$ are distinct modulo $p^2$, an application of Hensel's lemma shows that, for each fixed 
choice of ${\boldsymbol y}$, there are at most $k!(p^k)^{k(k-1)/2}\cdot p^{k(k-1)/2}$ choices for 
${\boldsymbol x}\mmod{p^{k^2}}$ satisfying (\ref{woolb.6}). Here, the factor $p^{k(k-1)/2}$ reflects 
the fact that, even though $x_i\not \equiv x_j\mmod{p^2}$ for $i\ne j$, one has 
$x_i\equiv x_j\mmod{p}$ for all $i$ and $j$. This situation is entirely analogous to that delivering 
(\ref{woola.6}) above, and thus we obtain
\[ K_{1,k}(X)\ll (M^{k+1})^{k(k-1)/2}\max_{\xi,\eta}\oint \biggl( \sum_{\substack{\nu=1\\ 
\nu\equiv \xi\mmod{p}}}^{p^{k^2}}|\grf_{k^2}(\bfalp;\nu)|^2\biggl)^k
|\grF_k(\bfalp;\eta)|^{2u}\diff \bfalp .\]

\par From here, as above, suitable applications of H\"older's inequality show that
\begin{equation}\label{woolb.6a}
\llbracket K_{1,k}(X)\rrbracket \ll \llbracket J_t(X/M^k)\rrbracket^{1-1/u} 
\llbracket K_{k,k^2}(X)\rrbracket^{1/u}.
\end{equation}
By substituting this estimate into (\ref{woolb.5}), we obtain the new upper bound
\[ \llbracket J_t(X)\rrbracket \ll \bigl( \llbracket J_t(X/M)\rrbracket \llbracket 
J_t(X/M^k)\rrbracket^{1/u}\bigr)^{1-1/u}\llbracket K_{k,k^2}(X)\rrbracket^{1/u^2}.\]
By iterating this process $N$ times, one obtains the relation
\begin{equation}\label{woolb.7}
\llbracket J_t(X)\rrbracket \ll \biggl( \prod_{r=0}^{N-1}\llbracket J_t(X/M^{k^r})\rrbracket^{1/u^r}
\biggr)^{1-1/u}\llbracket K_{k^{N-1},k^N}(X)\rrbracket^{1/u^N}.
\end{equation}
While this is a vaste oversimplification of what is actually established, it correctly identifies the relationship 
which underpins the efficient congruencing method.\par

Since $M^{k^N}<X^{1/3}$, our earlier discussion ensures that
\[ \llbracket J_t(X)\rrbracket \gg X^{\Lambda-\delta}\quad \text{and}\quad 
\llbracket J_t(X/M^{k^r})\rrbracket \ll (X/M^{k^r})^{\Lambda+\delta}\quad (0\le r\le N).\]
Meanwhile, an application of H\"older's inequality provides the trivial bound
\[ \llbracket K_{k^{N-1},k^N}(X)\rrbracket \ll \bigl(M^{k(k+1)/2}\bigr)^{k^N}X^{\Lambda+\delta}. \]
By substituting these estimates into (\ref{woolb.7}), we deduce that
\[ X^{\Lambda-\delta}\ll \biggl( X^{1/u^N}\prod_{r=0}^{N-1}
\bigl(X/M^{k^r}\bigr)^{(1-1/u)/u^r}\biggr)^{\Lambda+\delta}
\left( M^{k(k+1)/2}\right)^{(k/u)^N},\]
and hence $X^{\Lambda-\delta}\ll X^{\Lambda+\delta }(M^\Theta)^{(k/u)^N}$, where
\[ \Theta=\tfrac{1}{2}k(k+1)-(1-1/u)(\Lambda+\delta)\sum_{r=1}^N(u/k)^r.\]
But we have $u\ge k$, and so our hypotheses concerning $N$ and $\delta$ ensure that
\[ \Theta \le \tfrac{1}{2}k(k+1)-N(1-1/u)(\Lambda+\delta)<-\tfrac{1}{2}N\Lambda<-3(u/k)^N
\delta/\theta .\]
We therefore conclude that $X^{\Lambda-\delta}\ll X^{\Lambda+\delta}M^{-3\delta/\theta}\ll 
X^{\Lambda-2\delta}$. This relation yields the contradiction that establishes the desired conclusion 
$\Lambda=0$. We may therefore conclude that whenever $t\ge k(k+1)$, one has 
$J_t(X)\ll X^{2t-k(k+1)/2+\eps}$.\par

We have sketched the proof of the Main Conjecture for $J_t(X)$ when $t\ge k(k+1)$. Theorem 
\ref{wooltheorem3.1}, which represents the latest state of play in the efficient congruencing method, goes 
considerably further. Two ideas underpin these advances.\par

First, one may sacrifice some of the power potentially available from systems of congruences such as 
(\ref{woola.5}) or (\ref{woolb.6}) in order that the efficient congruencing method be applicable when 
$t<k(k+1)$. Let $r$ be a parameter with $2\le r\le k$, and define the generating function 
$\grF_c^{(r)}(\bfalp;\xi)$ by analogy with $\grF_c(\bfalp;\xi)$, though with $r$ (in place of $k$) 
underlying exponential sums $\grf_{c+1}(\bfalp;\xi_i)$. One may imitate the basic argument sketched 
above, with $t=(u+1)r$, to bound the analogue $K_{a,b}^{(r)}(X)$ of the mean value $K_{a,b}(X)$. In 
place of (\ref{woolb.6}) one now obtains the congruences
\begin{equation}\label{woolb.8}
(x_1-\eta)^j+\ldots +(x_r-\eta)^j\equiv (y_1-\eta)^j+\ldots +(y_r-\eta)^j\mmod{p^{jb}}
\quad (1\le j\le k).
\end{equation}
For simplicity, suppose that $r\le (k-1)/2$. Then by considering the $r$ congruence relations of highest 
degree here, one finds from Hensel's lemma that, for each fixed choice of ${\boldsymbol y}$, there are at 
most $k!$ choices for ${\boldsymbol x}\mmod{p^{(k-r)b}}$ satisfying (\ref{woolb.8}). Although this is a 
weaker congruence constraint than before on ${\boldsymbol x}$ and ${\boldsymbol y}$, the cost in terms 
of the number of choices is smaller, and so useful estimates may nonetheless be obtained for $J_t(X)$. 
Ideas along these lines underpin both the work \cite{Woo2013} of the author, and in the sharper form 
sketched above, that of Ford and the author \cite{FW2014}.\par

The second idea conveys us to the threshold of the Main Conjecture. Again we consider the mean values 
$K_{a,b}^{(r)}(X)$, and for simplicity put $r=k-1$. The congruences (\ref{woolb.8}) yield a constraint on 
the variables tantamount to $x_i\equiv y_i\mmod{p^{2b}}$ at little cost. Encoding this constraint using 
exponential sums, and applying H\"older's inequality, one bounds $K_{a,b}^{(k-1)}(X)$ in terms of 
$K_{a,b}^{(k-2)}(X)$ and $K_{b,2b}^{(k-1)}(X)$. Iterating this process to successively estimate 
$K_{a,b}^{(k-j)}(X)$ for $j=1,2,\ldots ,k-1$, we obtain a bound for $K_{a,b}^{(k-1)}(X)$ in terms of 
$K_{b,jb}^{(k-1)}(X)$ $(2\le j\le k)$ and $J_t(X/M^b)$. The heuristic potential of this idea amounts 
to a relation of the shape
\begin{equation}\label{woolb.9}
\llbracket K_{a,b}^{(k-1)}(X)\rrbracket \ll \biggl( \prod_{j=2}^k
\llbracket K_{b,jb}^{(k-1)}(X)\rrbracket^{\phi_j}\biggr)\llbracket J_t(X/M^b)\rrbracket^{1-(k-1)/s},
\end{equation}
where the exponents $\phi_j$ are approximately equal to $1/s$. Again, this substantially oversimplifies the 
situation, since non-negligible additional factors occur. However, one discerns a critical advantage over 
earlier relations such as (\ref{woolb.6a}). As one iterates (\ref{woolb.9}), one bounds 
$\llbracket K_{a,b}^{(k-1)}(X)\rrbracket $ in terms of new expressions 
$\llbracket K_{b,b'}^{(k-1)}(X)\rrbracket$, where the ratio $b'/b$ is on average about $\tfrac{1}{2}k+1$, 
as opposed to the previous ratio $k$. The relation (\ref{woolb.9}) may be converted into a substitute for 
(\ref{woolb.7}) of the shape
\[ \llbracket J_t(X)\rrbracket \ll \biggl( \prod_{r=0}^{N-1}\llbracket J_t(X/M^{\rho^r})
\rrbracket^{1/u^r}\biggr)^{1-1/u}\llbracket K^{(k-1)}_{\rho^{N-1},\rho^N}(X)\rrbracket^{1/u^N},\]
in which $\rho$ is close to $\tfrac{1}{2}k+1$ and $t=(u+1)(k-1)$. Thus, when $u\ge \rho$, we find as 
before that the lower bound $\llbracket J_t(X)\rrbracket \gg X^{\Lambda -\delta}$ is tenable only when 
$\Lambda=0$, and we have heuristically established the Main Conjecture when $t$ is only 
slightly larger than $k(k+1)/2$. Of course, the relation (\ref{woolb.9}) represents an idealised situation, 
and the proof in detail of the results in \cite{Woo2014a, Woo2014b} contains numerous complications 
requiring the resolution of considerable technical difficulties.

\section{Waring's problem} Investigations concerning the validity of the anticipated asymptotic formula in 
Waring's problem have historically followed one of two paths, associated on the one hand with Weyl, and 
on the other with Vinogradov. We recall our earlier notation, writing $R_{s,k}(n)$ for the number of 
representations of the natural number $n$ in the shape $n=x_1^k+\ldots +x_s^k$, with 
${\boldsymbol x}\in \dbN^s$. A heuristic application of the circle method suggests that for $k\ge 3$ and 
$s\ge k+1$, one should have
\begin{equation}\label{wool4.1}
R_{s,k}(n)=\frac{\Gamma (1+1/k)^s}{\Gamma (s/k)}\grS_{s,k}(n)n^{s/k-1}+o(n^{s/k-1}),
\end{equation}
where
\[ \grS_{s,k}(n)=\sum_{q=1}^\infty \sum^q_{\substack{a=1\\ (a,q)=1}}
\biggl( q^{-1}\sum_{r=1}^qe(ar^k/q)\biggr)^se(-na/q).\]
Under modest congruence conditions, one has $1\ll \grS_{s,k}(n)\ll n^\eps$, and thus the conjectural 
relation (\ref{wool4.1}) may be seen as an honest asymptotic formula (see 
\cite[\S\S4.3, 4.5 and 4.6]{Vau1997} for details). Let $\Gtil(k)$ denote the least integer $t$ with the 
property that, whenever $s\ge t$, the asymptotic formula (\ref{wool4.1}) holds for all large enough $n$.

\par Leaving aside the smallest exponents $k=1$ and $2$ accessible to classical methods, the first to 
obtain a bound for $\Gtil(k)$ were Hardy and Littlewood \cite{HL1922}, who devised a method based on 
Weyl differencing to show that $\Gtil(k)\le (k-2)2^{k-1}+5$. In 1938, Hua \cite{Hua1938} obtained a 
refinement based on the estimate
\begin{equation}\label{wool4.2}
\int_0^1|g_k(\alpha;X)|^{2^k}\diff \alpha \ll X^{2^k-k+\eps},
\end{equation}
in which $g_k(\alpha;X)$ is defined via (\ref{wool1.6}), showing that $\Gtil(k)\le 2^k+1$. For small 
values of $k$, this estimate remained the strongest known for nearly half a century. Finally, Vaughan 
\cite{Vau1986a, Vau1986b} succeeded in wielding Hooley's $\Delta$-functions to deduce that 
$\Gtil(k)\le 2^k$ for $k\ge 3$. For slightly larger exponents $k\ge 6$, this bound was improved by 
Heath-Brown \cite{HB1988} by combining Weyl differencing with a novel cubic mean value estimate. His 
bound $\Gtil(k)\le \frac{7}{8}2^k+1$ was, in turn, refined by Boklan \cite{Bok1994}, who exploited 
Hooley's $\Delta$-functions in this new setting to deduce that $\Gtil(k)\le \frac{7}{8}2^k$ for $k\ge 6$.

\par Turning now to large values of $k$, the story begins with Vinogradov \cite{Vin1935}, who showed 
that $\Gtil(k)\le 183k^9(\log k+1)^2$, reducing estimates previously exponential in $k$ to polynomial 
bounds. As Vinogradov's mean value theorem progressed to the state essentially captured by Theorem 
\ref{wooltheorem1.1}, bounds were rapidly refined to the form $\Gtil(k)\le (C+o(1))k^2\log k$, 
culminating in 1949 with Hua's bound \cite{Hua1949} of this shape with $C=4$. The connection with 
Vinogradov's mean value theorem is simple to explain, for on considering the underlying Diophantine 
systems, one finds that
\[ \int_0^1|g_k(\alpha;X)|^{2s}\diff \alpha =\sum_{\boldsymbol h}\oint |f_k(\bfalp;X)|^{2s}
e(-h_1\alp_1-\ldots -h_{k-1}\alp_{k-1})\diff \bfalp ,\]
where the summation is over $|h_j|\le sX^j$ $(1\le j\le k-1)$. The bound (\ref{wool1.4}) therefore leads 
via the triangle inequality and (\ref{wool1.2}) to the estimate
\begin{equation}\label{wool4.3}
\int_0^1|g_k(\alpha;X)|^{2s}\diff \alpha \ll X^{k(k-1)/2}J_{s,k}(X)\ll X^{2s-k+\Delta_{s,k}},
\end{equation}
which serves as a surrogate for (\ref{wool4.2}). In 1992, the author reduced the permissible value of $C$ 
from $4$ to $2$ by applying the repeated efficient differencing method \cite{Woo1992b}. A more efficient 
means of utilising Vinogradov's mean value theorem to bound $\Gtil(k)$ was found by Ford 
\cite{For1995} (see also \cite{Ust1998}), showing that $C=1$ is permissible. Refinements for smaller 
values of $k$ show that this circle of ideas surpasses the above-cited bound $\Gtil(k)\le \frac{7}{8}2^k$ 
when $k\ge 9$ (see Boklan and Wooley \cite{BW2012}).\par

We summarise the classical state of affairs in the following theorem.

\begin{theorem}[Classical status of $\Gtil(k)$]\label{wooltheorem4.1}
One has:
\begin{enumerate}
\item[(i)] $\Gtil(k)\le 2^k$ $(k=3,4,5)$ and $\Gtil(k)\le \tfrac{7}{8}2^k$ $(k=6,7,8)$;
\item[(ii)] $\Gtil(9)\le 365$, $\Gtil(10)\le 497$, $\Gtil(11)\le 627$, $\Gtil(12)\le 771$, ...;
\item[(iii)] $\Gtil(k)\le (1+o(1))k^2\log k$ $(k\ \text{large})$.
\end{enumerate}
\end{theorem}

The most immediate impact of the new efficient congruencing method in Vinogradov's mean value 
theorem \cite{Woo2012a} was the bound $\Gtil(k)\le 2k^2+2k-3$, valid for $k\ge 2$. This already 
supersedes the previous work presented in Theorem \ref{wooltheorem4.1} when $k\ge 7$. In particular, 
the obstinate factor of $\log k$ is definitively removed for large values of $k$. Subsequent refinements 
\cite{FW2014, Woo2012b, Woo2013, Woo2014a, Woo2014b} have delivered further progress, 
especially for smaller values of $k$, which we summarise as follows.

\begin{theorem}[Status of $\Gtil(k)$ after efficient congruencing]\label{wooltheorem4.2} One has:
\begin{enumerate}
\item[(i)] $\Gtil(k)\le 2^k$ $(k=3,4)$;
\item[(ii)] $\Gtil(5)\le 28$, $\Gtil(6)\le 43$, $\Gtil(7)\le 61$, $\Gtil(8)\le 83$, 
$\Gtil(9)\le 107$, $\Gtil(10)\le 134$, $\Gtil(11)\le 165$, $\Gtil(12)\le 199$, ...;
\item[(iii)] $\Gtil(k)\le (C+o(1))k^2$ $(k\ \text{large})$, where $C=1.54079$ is an approximation to the 
number $(5+6\xi-3\xi^2)/(2+6\xi)$, in which $\xi$ is the real root of $6\xi^3+3\xi^2-1$.
\end{enumerate}
\end{theorem}

A comparison of Theorems \ref{wooltheorem4.1} and \ref{wooltheorem4.2} reveals that the classical 
Weyl-based bounds have now been superseded for $k\ge 5$. The latest developments 
\cite{Woo2014a, Woo2014c} hint, indeed, at further progress even when $k=4$. These advances for 
smaller values of $k$ stem in part, of course, from the substantial progress in our new bounds for 
$J_{s,k}(X)$, as outlined in Theorems \ref{wooltheorem3.1} and \ref{wooltheorem3.2}. However, an 
important role is also played by a novel mean value estimate for moments of $g_k(\alp;X)$. Define the 
minor arcs $\grm=\grm_k$ to be the set of real numbers $\alp\in [0,1)$ satisfying the property that, 
whenever $a\in \dbZ$ and $q\in \dbN$ satisfy $(a,q)=1$ and $|q\alp-a|\le X^{1-k}$, then $q>X$. The 
argument of the proof of \cite[Theorem 2.1]{Woo2012b} yields the bound
\begin{equation}\label{wool4.4}
\int_\grm |g_k(\alp;X)|^{2s}\diff \alp \ll X^{\frac{1}{2}k(k-1)-1}(\log X)^{2s+1}J_{s,k}(X).
\end{equation}
We thus infer from Theorem \ref{wooltheorem3.1}(iii) that whenever $k\ge 3$ and $s\ge k(k-1)$, then
\[ \int_\grm |g_k(\alp;X)|^{2s}\diff \alp \ll X^{2s-k-1+\eps}.\]
As compared to the classical approach embodied in (\ref{wool4.3}), an additional factor $X$ has been 
saved in these estimates at no cost in terms of the number of variables, and for smaller values of $k$ this 
is a very substantial gain.\par

For large values of $k$, the enhancement of Ford \cite{For1995} given by Ford and Wooley 
\cite[Theorem 8.5]{FW2014} remains of value. When $k,s\in \dbN$, denote by $\eta(s,k)$ the least 
number $\eta$ with the property that, whenever $X$ is sufficiently large 
in terms of $s$ and $k$, one has $J_{s,k}(X)\ll_\eps X^{2s-k(k+1)/2+\eta+\eps}$. Let $r\in \dbN$ 
satisfy $1\le r\le k-1$. Then \cite[Theorem 8.5]{FW2014} shows that whenever $s\ge r(r-1)/2$, one has
\[ \int_0^1 |g_k(\alp;X)|^{2s}\diff \alp \ll X^{2s-k+\eps}\bigl( X^{\eta_r^*(s,k)-1/r}+
X^{\eta_r^*(s,k-1)}\bigr) ,\]
where $\eta_r^*(s,w)=r^{-1}\eta(s-r(r-1)/2,w)$ for $w=k-1,k$.\par

Finally, we note that familiar conjectures concerning mean values of the exponential sum $g_k(\alp;X)$ 
imply that one should have $\Gtil(k)\le 2k+1$ for each $k\ge 3$, and indeed it may even be the case that 
$\Gtil(k)\le k+1$.

\section{Estimates of Weyl-type, and distribution mod $1$} Pointwise estimates for exponential sums 
appear already in the work of Weyl \cite{Wey1916} in 1916. By applying $k-1$ Weyl-differencing steps, 
one bounds the exponential sum $f_k(\bfalp;X)$ in terms of a new exponential sum over a linear 
polynomial, and this may be estimated by summing what is, after all, a geometric progression. In this way, 
one obtains the classical version of Weyl's inequality (see \cite[Lemma 2.4]{Vau1997}).

\begin{theorem}[Weyl's inequality]\label{wooltheorem5.1}
Let $\bfalp\in \dbR^k$, and suppose that $a\in \dbZ$ and $q\in \dbN$ satisfy $(a,q)=1$ and 
$|\alp_k-a/q|\le q^{-2}$. Then one has
\begin{equation}\label{wool5.1}
|f_k(\bfalp;X)|\ll X^{1+\eps}\left(q^{-1}+X^{-1}+qX^{-k}\right)^{2^{1-k}}.
\end{equation}
\end{theorem}

This provides a non-trivial estimate for $f_k(\bfalp;X)$ when the leading coefficient $\alp_k$ is not 
well-approximated by rational numbers. Consider, for example, the set $\grm=\grm_k$ defined in the 
preamble to (\ref{wool4.4}). When $\alp_k\in \grm$, an application of Dirichlet's theorem on Diophantine 
approximation shows that there exist $a\in \dbZ$ and $q\in \dbN$ with $(a,q)=1$ such that 
$q\le X^{k-1}$ and $|q\alp-a|\le X^{1-k}$. The definition of $\grm$ then implies that $q>X$, and so 
Theorem \ref{wooltheorem5.1} delivers the bound
\begin{equation}\label{wool5.2}
\sup_{\alp_k\in \grm}|f_k(\bfalp;X)|\ll X^{1-\sigma(k)+\eps},
\end{equation}
in which $\sigma(k)=2^{1-k}$. Heath-Brown's variant \cite[Theorem 1]{HB1988} of Weyl's inequality 
applies mean value estimates for certain cubic exponential sums that, for $k\ge 6$, give bounds superior 
to (\ref{wool5.1}) when $q$ lies in the range $X^{5/2}<q<X^{k-5/2}$. By making use of the cubic case 
of the Main Conjecture in Vinogradov's mean value theorem \cite{Woo2014c}, the author \cite{Woo2014d} 
has extended this range to $X^2<q<X^{k-2}$.

\begin{theorem}\label{wooltheorem5.2}
Let $k\ge 6$, and suppose that $\alp\in \dbR$, $a\in \dbZ$ and $q\in \dbN$ satisfy $(a,q)=1$ and 
$|\alp-a/q|\le q^{-2}$. Then one has
\[ |g_k(\alp;X)|\ll X^{1+\eps}\Theta^{2^{-k}}+X^{1+\eps}(\Theta /X)^{\frac{2}{3}2^{-k}},\]
where $\Theta =q^{-1}+X^{-3}+qX^{-k}$.
\end{theorem}

For comparison, we note that Heath-Brown \cite[Theorem 1]{HB1988} obtains the bound
$|g_k(\alp;X)|\ll X^{1+\eps}(X\Theta)^{\frac{4}{3}2^{-k}}$. Robert and Sargos 
\cite[Th\'eor\`eme 4 et Lemme 7]{RS2000} extend these ideas when $k\ge 8$ to show that 
$|g_k(\alp;X)|\ll X^{1+\eps} \bigl( X^{17/8}\Theta'\bigr)^{\frac{8}{5}2^{-k}}$, in which 
$\Theta'=q^{-1}+X^{-4}+qX^{-k}$. See Parsell \cite{Par2014} for a refinement when $k=8$.\par

The above methods yield exponents exponentially small in $k$. By substituting estimates for $J_{s,k}(X)$ 
into the conclusion of Theorem \ref{wooltheorem1.3}, one obtains analogous bounds polynomial in $k$. 
Classical versions of Vinogradov's mean value theorem yield estimates of the shape (\ref{wool5.2}) with 
$\sigma(k)^{-1}=(C+o(1))k^2\log k$. Thus, Linnik \cite{Lin1943} obtained the permissible value 
$C=22\,400$, and Hua \cite{Hua1949} obtained $C=4$ in 1949. This was improved via efficient 
differencing \cite{Woo1992b} in 1992 to $C=2$, and subsequently the author \cite{Woo1995} obtained 
$C=3/2$ by incorporating some ideas of Bombieri \cite{Bom1990}. The latest developments in efficient 
congruencing yield the new exponent $\sigma(k)^{-1}=2(k-1)(k-2)$ for $k\ge 3$, a conclusion that 
removes the factor $\log k$ from earlier estimates, and improves on Weyl's inequality for $k\ge 7$.

\begin{theorem}\label{wooltheorem5.3}
Let $k$ be an integer with $k\ge 3$, and let $\bfalp\in \dbR^k$. Suppose that there exists a natural 
number $j$ with $2\le j\le k$ such that, for some $a\in \dbZ$ and $q\in \dbN$ with $(a,q)=1$, one has 
$|\alp_j-a/q|\le q^{-2}$. Then one has
\[ |f_k(\bfalp;X)|\ll X^{1+\eps} (q^{-1}+X^{-1}+qX^{-j})^{\sigma(k)},\]
where $\sigma(k)^{-1}=2(k-1)(k-2)$.
\end{theorem}

This conclusion makes use of Theorem \ref{wooltheorem3.1}(iii), and improves slightly on 
\cite[Theorem 11.1]{Woo2014a}. When $k\ge 6$ and $\alp$ lies on a suitable subset of $\dbR$, the 
above-cited work of Heath-Brown may deliver estimates for $|g_k(\alp;X)|$ superior to those stemming 
from Weyl's inequality, and similar comments apply to the work of Robert and Sargos, and of Parsell, when 
$k\ge 8$. However, Theorem \ref{wooltheorem5.3} proves superior in all circumstances to the estimates 
of Heath-Brown when $k>7$, and to the estimates of Robert and Sargos for all exponents $k$.\par

In many applications, it is desirable to have available estimates for $|f_k(\bfalp;X)|$ that depend on 
simultaneous approximations to $\alp_1,\ldots ,\alp_k$ of a given height. This is a subject to which 
R. C. Baker and W. M. Schmidt have made significant contributions. By exploiting such methods in 
combination with our new estimates for $J_{s,k}(X)$, one obtains the following conclusion (compare 
\cite[Theorem 11.2]{Woo2014a}).

\begin{theorem}\label{wooltheorem5.4}
Let $k$ be an integer with $k\ge 3$, and let $\tau$ and $\delta$ be real numbers with 
$\tau^{-1}>4(k-1)(k-2)$ and $\delta>k\tau$. Suppose that $X$ is sufficiently large in terms of $k$, 
$\delta$ and $\tau$, and further that $|f_k(\bfalp;X)|>X^{1-\tau}$. Then there exist integers 
$q,a_1,\ldots ,a_k$ such that $1\le q\le X^\delta$ and $|q\alp_j-a_j|\le X^{\delta -j}$ $(1\le j\le k)$.
\end{theorem}

Here, the constraint $\tau^{-1}>4(k-1)(k-2)$ may be compared with the corresponding hypothesis 
$\tau^{-1}>(8+o(1))k^2\log k$ to be found in \cite[Theorem 4.5]{Bak1986}, and the Weyl-based 
bound $\tau^{-1}>2^{k-1}$ obtained in \cite[Theorem 5.2]{Bak1986} (see also \cite{Bak1992}). The 
conclusion of Theorem \ref{wooltheorem5.4} is superior to the latter for $k>8$.\par

Bounds for exponential sums of Weyl-type may be converted into equidistribution results for polynomials 
modulo $1$ by applying estimates of Erd\H os-Tur\'an type. Write $\|\theta\|$ for the least value of 
$|\theta -n|$ for $n\in \dbZ$, and consider a sequence $(x_n)_{n=1}^\infty$ of real numbers . Then it 
follows from \cite[Theorem 2.2]{Bak1986}, for example, that whenever $\|x_n\|\ge M^{-1}$ for 
$1\le n\le N$, then
\[ \sum_{1\le m\le M}\biggl| \sum_{n=1}^Ne(mx_n)\biggr| >\tfrac{1}{6}N.\]
By carefully exploiting this result using the methods of Baker \cite{Bak1986}, one deduces from Theorem 
\ref{wooltheorem5.4} the following conclusion (compare \cite[Theorem 11.3]{Woo2014a}).

\begin{theorem}\label{wooltheorem5.5}
When $k\ge 3$, put $\tau(k)=1/\left( 4(k-1)(k-2)\right)$. Then whenever $\bfalp\in \dbR^k$ and $N$ is 
sufficiently large in terms of $k$ and $\eps$, one has
\[ \min_{1\le n\le N}\| \alp_1n+\alp_2n^2+\ldots +\alp_kn^k\|<N^{\eps-\tau(k)}.\]
\end{theorem}

\section{Further applications} Vinogradov's mean value theorem finds application in numerous 
number-theoretic problems, besides those discussed in the previous two sections. We take the opportunity 
now to outline several applications, emphasising recent developments.\vskip.1cm

\noindent \emph{(i) Tarry's problem.} When $h$, $k$ and $s$ are positive integers with $h\ge 2$, 
consider the Diophantine system
\begin{equation}\label{wool6.1}
\sum_{i=1}^sx_{i1}^j=\sum_{i=1}^sx_{i2}^j=\ldots =\sum_{i=1}^sx_{ih}^j\quad (1\le j\le k).
\end{equation}
Let $W(k,h)$ denote the least natural number $s$ having the property that the simultaneous equations 
(\ref{wool6.1}) possess an integral solution ${\boldsymbol x}$ with    
\[ \sum_{i=1}^sx_{iu}^{k+1}\ne \sum_{i=1}^sx_{iv}^{k+1}\quad (1\le u<v\le h).\]
The problem of estimating $W(k,h)$ was intensely investigated by E. M. Wright and L.-K. Hua (see 
\cite{Hua1938b, Hua1949b, Wri1948}), the latter obtaining $W(k,h)\le k^2(\log k+O(1))$ for $h\ge 2$. 
The argument of the proof of \cite[Theorem 1.3]{Woo2012a} shows that $W(k,h)\le s$ whenever one can 
establish the estimate $J_{s,k+1}(X)=o(X^{2s-k(k+1)/2})$. By using this criterion together with the 
estimates for $J_{s,k+1}(X)$ obtained via the latest efficient congruencing methods, one obtains 
substantial improvements in these earlier conclusions (see \cite[Theorem 12.1]{Woo2014a} and 
\cite[Theorem 12.1]{Woo2014b}).

\begin{theorem}\label{wooltheorem6.1}
When $h$ and $k$ are natural numbers with $h\ge 2$ and $k\ge 3$, one has 
$W(k,h)\le \tfrac{5}{8}(k+1)^2$. Moreover, when $k$ is large, one has 
$W(k,h)\le \tfrac{1}{2}k(k+1)+1$.
\end{theorem}

Although the last of these conclusions achieves the limit of current analytic approaches to bounding 
$W(k,h)$, explicit numerical examples are available\footnote{See the website 
http://euler.free.fr/eslp/eslp.htm.} which may be applied to show that $W(k,2)=k+1$ for $1\le k\le 9$ 
and $k=11$.\vskip.1cm

\noindent \emph{(ii) Sum-product theorems.} When $A$ is a finite set of real numbers, define the sets 
$A+A=\{ x+y:x,y\in A\}$ and $A\cdot A=\{ xy:x,y\in A\}$, and more generally
\[ hA=\{x_1+\ldots +x_h:{\boldsymbol x}\in A^h\} \quad \text{and}\quad 
A^{(h)}=\{ x_1\ldots x_h:{\boldsymbol x}\in A^h\} .\]
A conjecture of Erd\H os and Szemer\'edi \cite{ES1983} asserts that for any finite set of integers $A$, one 
has $|A+A|+|A\cdot A|\gg_\eps |A|^{2-\eps}$. It is also conjectured that whenever $A$ is a finite set of 
\textit{real} numbers, then for each $h\in \dbN$, one should have 
$|hA|+|A^{(h)}|\gg_{\eps,h}|A|^{h-\eps}$. Chang \cite{Cha2003} has made progress towards 
this conjecture by showing that when $A$ is a finite set of integers, and $|A.A|<|A|^{1+\eps}$, 
then $|hA|\gg_{\eps,h} |A|^{h-\delta}$, where $\delta\rightarrow 0$ as $\eps\rightarrow 0$. 
Subsequently, Bourgain and Chang \cite{BC2003} showed that for any $b\ge 1$, there exists $h\ge 1$ 
with the property that $|hA|+|A^{(h)}|\gg |A|^b$. By exploiting bounds for $W(k,h)$ of the type 
given by Theorem \ref{wooltheorem6.1}, Croot and Hart \cite{CH2010} have made progress toward an 
analogue of such conclusions for sets of real numbers.

\begin{theorem}\label{wooltheorem6.2}
Suppose that $\eps>0$ and $|A\cdot A|\le |A|^{1+\eps}$. Then there exists a number $\lambda>0$ such that, when $h$ is large enough in terms of $\eps$, one has 
$|h(A\cdot A)|>|A|^{\lambda h^{1/3}}$.
\end{theorem}

This conclusion (see \cite[Theorem 11.5]{Woo2013}) improves on \cite[Theorem 2]{CH2010}, where a 
similar result is obtained with $(h/\log h)^{1/3}$ in place of the exponent $h^{1/3}$.
\vskip.1cm

\noindent \emph{(iii) The Hilbert-Kamke problem and its brethren.} Hilbert \cite{Hil1909} considered an 
extension of Waring's problem related to Vinogradov's mean value theorem. When 
$n_1,\ldots ,n_k\in \dbN$, let $R_{s,k}({\boldsymbol n})$ denote the number of solutions of the system
\begin{equation}\label{wool6.2}
x_1^j+\ldots +x_s^j=n_j\quad (1\le j\le k),
\end{equation}
with ${\boldsymbol x}\in \dbN^s$. Put $\displaystyle{X=\max_{1\le j\le k}n_j^{1/j}}$, and then write
\[ {\cal J}_{s,k}({\boldsymbol n})=\int_{\dbR^k}\biggl( 
\int_0^1e(\beta_1\gamma +\ldots +\beta_k\gamma^k)\diff \gamma \biggr)^s
e(-\beta_1n_1/X-\ldots -\beta_kn_k/X^k) \diff {\boldsymbol \beta} \]
and
\[ {\cal S}_{s,k}({\boldsymbol n})=\sum_{q=1}^\infty 
\sum_{\substack{1\le a_1,\ldots ,a_k\le q\\ (q,a_1,\ldots ,a_k)=1}}
\left( q^{-1}f_k({\boldsymbol a}/q,q) \right)^se\left( -(a_1n_1+\ldots +a_kn_k)/q\right) .\]
See \cite{Ark1984, Mit1986, Mit1987} for an account of the analysis of this problem, and in particular for a 
discussion of the conditions under which real and $p$-adic solutions exist for the system (\ref{wool6.2}). 
While the conditions $n_k^{j/k}\le n_j\le s^{1-j/k}n_k^{j/k}$ $(1\le j\le k)$ are plainly necessary, one 
finds that $p$-adic solubility is not assured when $s<2^k$. This classical technology gives an asymptotic 
formula for $R_{s,k}({\boldsymbol n})$ provided that $s\ge (4+o(1))k^2\log k$. Efficient congruencing 
methods lead to considerable progress. The following result improves on \cite[Theorem 9.2]{Woo2012a} 
using Theorem \ref{wooltheorem3.1}(iii).

\begin{theorem}\label{wooltheorem6.3}
Let $s,k\in \dbN$ and ${\boldsymbol n}\in \dbN^k$. Suppose that $X=\max n_j^{1/j}$ is sufficiently 
large in terms of $s$ and $k$, and that the system (\ref{wool6.2}) has non-singular real and $p$-adic 
solutions. Then whenever $k\ge 3$ and $s\ge 2k^2-2k+1$, one has
\[ R_{s,k}({\boldsymbol n})={\cal J}_{s,k}({\boldsymbol n}){\cal S}_{s,k}({\boldsymbol n})
X^{s-k(k+1)/2}+o(X^{s-k(k+1)/2}).\]
\end{theorem}

Similar arguments apply to more general Diophantine systems. Let $k_1,\ldots ,k_t$ be distinct positive 
integers. Suppose that $s, k\in \dbN$, and that $a_{ij}\in \dbZ$ for $1\le i\le t$ and $1\le j\le s$. Write
\[ \phi_i({\boldsymbol x})=a_{i1}x_1^{k_i}+\ldots +a_{is}x_s^{k_i}\quad (1\le i\le t),\]
and consider the Diophantine system $\phi_i({\boldsymbol x})=0$ $(1\le i\le t)$. We write 
$N(B;{\boldsymbol \phi})$ for the number of integral solutions of this system with 
$|{\boldsymbol x}|\le B$. When $L>0$, define
\[ \sigma_\infty =\lim_{L\rightarrow \infty} \int_{|{\boldsymbol \xi}|\le 1}\prod_{i=1}^t
\max \left\{ 0, L(1-L|\phi_i({\boldsymbol \xi})|)\right\}\diff {\boldsymbol \xi}.\]
Also, for each prime number $p$, put
\[ \sigma_p=\lim_{H\rightarrow \infty} p^{H(t-s)}\text{card}\{ {\boldsymbol x}\in (\dbZ/p^H\dbZ)^s: 
\phi_i({\boldsymbol x})\equiv 0\ (\text{mod}\, p^H)\ (1\le i\le t)\} .\]
By applying the Hardy-Littlewood method, a fairly routine application of Theorem 
\ref{wooltheorem3.1}(iii) delivers the following conclusion (compare \cite[Theorem 9.1]{Woo2012a}).

\begin{theorem}\label{wooltheorem6.4}
Let $s$ and $k$ be natural numbers with $k\ge 3$ and $s\ge 2k^2-2k+1$. Suppose that $\max k_i\le k$, 
and that $a_{ij}$ $(1\le i\le t,\, 1\le j\le s)$ are non-zero integers. Suppose in addition that the system of 
equations $\phi_i({\boldsymbol x})=0$ $(1\le i\le t)$ has non-singular real and $p$-adic solutions, for 
each prime number $p$. Then
\[ N(B;{\boldsymbol \phi})\sim \sigma_\infty \biggl( \prod_p \sigma_p\biggr) B^{s-k_1-\ldots -k_t}.\]
\end{theorem}
\vskip.1cm

\noindent \emph{(iv) Solutions of polynomial congruences in short intervals.} There has been much 
activity in recent years concerning the solubility of polynomial congruences in short intervals, some of 
which makes use of estimates associated with Vinogradov's mean value theorem. Let 
$f\in {\mathbb F}_p[X]$ have degree $m\ge 3$, and let $M$ be a positive integer with $M<p$. Denote by 
$I_f(M;R,S)$ the number of solutions of the congruence $y^2\equiv f(x)\mmod{p}$, with 
$(x,y)\in [R+1,R+M]\times [S+1,S+M]$. Weil's bounds for exponential sums yield the estimate 
$I_f(M;R,S)=M^2p^{-1}+O(p^{1/2}(\log p)^2)$, one that is worse than trivial for 
$M\le p^{1/2}(\log p)^2$. The work of Chang et al. \cite{CCG2014} gives estimates that remain 
non-trivial for significantly smaller values of $M$.

\begin{theorem}\label{wooltheorem6.5}
Let $f\in {\mathbb F}_p[X]$ be any polynomial of degree $m\ge 4$. Then whenever $M$ is a positive 
integer with $1\le M<p$, we have
\[ I_f(M;R,S)\ll M^{1+\eps}\left( M^3p^{-1}+M^{3-m}\right)^{1/(2k(k-1))}.\]
\end{theorem}

This follows from \cite[Theorem 4]{CCG2014} on applying Theorem \ref{wooltheorem3.1}(iii). In 
particular, for any $\eps>0$, one finds that there exists a $\delta>0$, depending only on $\eps$ and 
deg$(f)$, such that whenever $M<p^{1/3-\eps}$ and deg$(f)\ge 4$, then $I_f(M;R,S)\ll M^{1-\delta}$.
\vskip.1cm

\noindent \emph{(v) The zero-free region for the Riemann zeta function.} We would be remiss not to 
mention the role of Vinogradov's mean value theorem in the proof of the widest available zero-free region 
for the Riemann zeta function. The sharpest estimates date from work of Vinogradov \cite{Vin1958} and 
Korobov \cite{Kor1958} in 1958 (see also \cite{Wal1963}). Thus, there is a positive constant $c_1$ with 
the property that $\zeta(s)\ne 0$ when $s=\sigma+it$, with $\sigma,t\in \dbR$, whenever $|t|\ge 3$ and 
$\sigma\ge 1-c_1(\log |t|)^{-2/3}(\log \log |t|)^{-1/3}$. More recently, Ford \cite{For2002} has shown 
that one may take $c_1=1/57.54$. This, in turn leads to an effective version of the prime number theorem 
of the shape
\[ \pi(x)=\int_2^x\frac{\diff t}{\log t}+O\left( x\exp\bigl(-c_2(\log x)^{3/5}(\log \log x)^{-1/5}\bigr) 
\right), \]
where $c_2=0.2098$. Using currently available methods, the nature of the constant $C(k,r)$ in estimates 
of the shape (\ref{wool1.4}) is significant for estimates of this type, while the precise nature of the defect 
in the exponent $\Delta_{s,k}$ less so. Thus, although the new estimates for $J_{s,k}(X)$ stemming from 
efficient congruencing have the potential to impact the numerial constants $c_1$ and $c_2$, the 
dependence on $t$ and $x$, respectively, in the above estimates has not been affected. 

\section{Generalisations} Thus far, we have focused on estimates for $J_{s,k}(X)$, the number of 
solutions (over the ring $\dbZ$) of the translation-dilation invariant system (\ref{wool1.3}) with 
$1\le {\boldsymbol x},{\boldsymbol y}\le X$. Previous authors have considered generalisations in which 
either the ring, or else the translation-dilation invariant system, is varied.\vskip.1cm

\noindent \emph{(i) Algebraic number fields.} The arguments underlying the proof of Theorem 
\ref{wooltheorem3.1} change little when the setting is shifted from $\dbZ$ to the ring of integers of a 
number field. When $s\ge k(k-1)$, the ensuing estimates are at most a factor $X^\eps$ away from the 
upper bound predicted by a heuristic application of the circle method. In common with Birch's use of Hua's 
lemma in number fields \cite{Bir1961}, our estimates are therefore robust to variation in the degree of the 
field extension, since Weyl-type estimates for exponential sums no longer play a significant role in 
applications. In forthcoming work we apply such ideas to establish the following result.

\begin{theorem}\label{wooltheorem9.1} Let $L/\dbQ$ be a field extension of finite degree. Suppose that 
$d\ge 3$, $s>2d(d-1)$ and ${\boldsymbol a}\in (L^\times)^s$. Then the hypersurface defined by 
$a_1x_1^d+\ldots +a_sx_s^d=0$ satisfies weak approximation 
and the Hasse principle over $L$.
\end{theorem}

For comparison, Birch \cite[Theorem 3]{Bir1961} gives such a conclusion only for $s>2^d$, 
while the work of K\"orner \cite{Kor1962} yields analogous conclusions in which the number of 
variables is larger, and depends also on the degree of the field extension $L/\dbQ$.\vskip.1cm

\noindent \emph{(ii) Function fields.} Consider a finite field $\dbF_q$ of characteristic $p$. Let 
$B\in \dbN$ be large enough in terms of $q$, $k$ and $s$, and denote by $J_{s,k}(B;q)$ the number of 
solutions of (\ref{wool1.3}) with $x_i,y_i\in \dbF_q[t]$ $(1\le i\le s)$ having degree at most $B$. When 
$p<k$, one can reduce (\ref{wool1.3}) to a minimal translation-invariant system in which certain 
equations are omitted. We write $K$ for the sum of the degrees of this minimal system, so that 
$K=k(k+1)/2$ when $p>k$, and $K<k(k+1)/2$ when $p<k$. Then, when $s\ge k(k+1)$, the efficient 
congruencing method adapts to give the upper bound $J_{s,k}(B;q)\ll (q^B)^{2s-K+\eps}$. This and 
much more is contained in forthcoming work of the author joint with Y.-R. Liu, generalisations of which are 
described in \cite{KLZ2014}.\vskip.1cm

\noindent \emph{(iii) Multidimensional analogues.} Vinogradov's methods have been generalised to 
multidimensional settings by Arkhipov, Chubarikov and Karatsuba \cite{ACK2004, AKC1980}, Parsell 
\cite{Par2005} and Prendiville \cite{Pre2013}. Variants of the efficient congruencing method deliver 
much sharper conclusions in far greater generality. Let $r,s,d\in \dbN$, and consider a linearly 
independent system of homogeneous polynomials ${\boldsymbol F}=(F_1,\ldots ,F_r)$, where 
$F_j({\boldsymbol z})\in \dbZ[z_1,\ldots ,z_d]$. Suppose that for $1\le j\le r$ and $1\le l\le d$, 
the polynomial $\partial F_j/\partial z_l$ lies in span$(1,F_1,\ldots ,F_r)$. Such a 
{\textit{reduced translation-dilation invariant}} system is said to have {\textit{rank}} $r$, 
{\textit{dimension}} $d$, {\textit{degree}} $k=\max \text{deg}\, F_j$, and {\textit{weight}} 
$K=\sum_1^r \text{deg}\, F_j$. Denote by $J_s(X;{\boldsymbol F})$ the number of integral solutions 
of the system of equations
\[ \sum_{i=1}^s\left( F_j({\boldsymbol x}_i)-F_j({\boldsymbol y}_i)\right)=0
\quad (1\le j\le r),\]
with $1\le {\boldsymbol x}_i,{\boldsymbol y}_i\le X$ $(1\le i\le s)$. The work of Parsell, Prendiville and 
the author \cite[Theorem 2.1]{PPW2013} provides a general estimate for $J_s(X;{\boldsymbol F})$ 
matching the predictions of the appropriate analogue of the Main Conjecture.

\begin{theorem}\label{wooltheorem9.2} Let ${\boldsymbol F}$ be a reduced translation-dilation invariant 
system of rank $r$, dimension $d$, degree $k$ and weight $K$. Then 
$J_s(X;{\boldsymbol F})\ll X^{2sd-K+\eps}$ for $s\ge r(k+1)$.
\end{theorem}

Reduced translation-dilation invariant systems are easy to generate by taking successive partial 
derivatives and reducing to a linearly independent spanning set. Thus, for example, the initial seed 
$x^5+3x^2y^3$ gives rise to just such a system
\[ {\boldsymbol F}=\{ x^5+3x^2y^3, 5x^4+6xy^3, 
x^2y^2,10x^3+3y^3,xy^2,x^2y,x^2,xy,y^2,x,y\},\]
with $d=2$, $r=11$, $k=5$, $K=30$. We therefore see from Theorem \ref{wooltheorem9.2} that 
$J_s(X;{\boldsymbol F})\ll X^{4s-30+\eps}$ for $s\ge 66$. Theorem \ref{wooltheorem9.2} should be 
susceptible to improvement by using the ideas underlying multigrade efficient congruencing 
\cite{Woo2014a, Woo2014b, Woo2014c}.

\section{Challenges} The remarkable success of the efficient congruencing method encourages 
ambitious speculation concerning other potential applications, a topic we briefly explore.\vskip.1cm

\noindent \emph{(i) The Main Conjecture for larger $s$.} In Theorem \ref{wooltheorem3.1}, one sees 
that the upper bound $J_{s,k}(X)\ll X^{s+\eps}$ predicted by the Main Conjecture is now known to hold 
for $1\le s\le \tfrac{1}{2}k(k+1)-t_k$, where $t_k=\tfrac{1}{3}k+O(k^{2/3})$. In striking contrast, on 
the other side of the critical value $s=\tfrac{1}{2}k(k+1)$, the upper bound 
$J_{s,k}(X)\ll X^{2s-k(k+1)/2+\eps}$ is known to hold only when $s\ge \tfrac{1}{2}k(k+1)+u_k$, 
where $u_k=\tfrac{1}{2}k(k-3)$. Plainly, the value of $u_k$ is substantially larger than $t_k$, and an 
intriguing possibility is that a hitherto unseen refinement of the method might reduce $u_k$ to a size 
more similar to that of $t_k$. This would have great significance in numerous applications.\vskip.1cm

\noindent \emph{(ii) Paucity.} When $k\ge 3$ and $1\le s<\tfrac{1}{2}k(k+1)$, we have precise 
asymptotics for $J_{s,k}(X)$ only when $s\le k+1$. Since the formula $J_{s,k}(X)=T_s(X)\sim s!X^s$ 
is trivial for $1\le s\le k$, the case $s=k+1$ is the only one with content. It is tempting to 
speculate that a suitable adaptation of efficient congruencing might confirm that 
$J_{s,k}(X)=T_s(X)+O(X^{s-\delta})$, for some $\delta>0$, for some exponent $s\ge k+2$.
\vskip.1cm

\noindent \emph{(iii) Minor arc bounds.} When $q\in \dbN$ and ${\boldsymbol a}\in \dbZ^k$, denote by 
$\grM(q,{\boldsymbol a})$ the set of points $\bfalp \in [0,1)^k$ such that $|q\alp_j-a_j|\le X^{1-j}$ 
$(1\le j\le k)$. Write $\grM$ for the union of the boxes $\grM(q,{\boldsymbol a})$ with 
$0\le a_j\le q\le X$ $(1\le j\le k)$ and $(q,a_1,\ldots ,a_k)=1$, and put $\grm=[0,1)^k\setminus \grM$. 
The methods of \S7 provide estimates of the shape $|f_k(\bfalp;X)|\ll X^{1-\sigma_k+\eps}$ for 
$\bfalp \in \grm$. However, when $s=k(k-1)+t$ and $t\ge 1$, our most efficient means of estimating 
moments of $f_k(\bfalp;X)$ of order $2s$, restricted to minor arcs, proceeds by applying Theorem 
\ref{wooltheorem3.1}(iii) via the trivial bound
\begin{align*}
\int_\grm|f_k(\bfalp;X)|^{2s}\diff \bfalp &\ll \Bigl( \sup_{\bfalp \in \grm}|f_k(\bfalp;X)|\Bigr)^{2t}
\oint |f_k(\bfalp;X)|^{2k(k-1)}\diff \bfalp \\
&\ll X^{2s-\frac{1}{2}k(k+1)-2t\sigma_k+\eps}.
\end{align*}
This bound is relatively weak, even when $t$ is large. Efficient congruencing provides a possible means of 
deriving estimates directly for such moments, and might even lead to improvements in our lower bounds 
for permissible exponents $\sigma_k$.\vskip.1cm

\noindent \emph{(iv) Non-translation invariant systems.} The system (\ref{wool1.3}) is translation-dilation 
invariant. A major desideratum is to apply a variant of efficient congruencing to systems of equations that 
are {\textit{not}} translation invariant. The author has forthcoming work applicable to systems that are 
only approximately translation invariant.



\begin{thebibliography}{17}

\bibitem{Ark1984}
Arkhipov, G. I., On the Hilbert-Kamke problem, \emph{Izv. Akad. Nauk SSSR Ser. Mat.} \textbf{48} 
(1984), no. 1, 3--52.

\bibitem{ACK2004}
Arkhipov, G. I., Chubarikov, V. N. and Karatsuba, A. A., \emph{Trigonometric sums in number theory and 
analysis}, Walter de Gruyter, Berlin, 2004.

\bibitem{AK1978}
Arkhipov, G. I. and Karatsuba, A. A., A new estimate of an integral of I. M. Vinogradov, 
\emph{Izv. Akad. Nauk SSSR Ser. Mat.} \textbf{42} (1978), no. 4, 751--762.

\bibitem{AKC1980}
Arkhipov, G. I., Karatsuba, A. A., and Chubarikov, V. N., Multiple Trigonometric Sums, \emph{Trudy Mat. 
Inst. Steklov} \textbf{151} (1980), 1--126. 

\bibitem{Bak1986}
Baker, R. C., \emph{Diophantine inequalites}, London Mathematical Society Monographs, vol. \textbf{1}, 
Oxford University Press, Oxford, 1986.

\bibitem{Bak1992}
Baker, R. C., Correction to: ``Weyl sums and Diophantine approximation'' [J. London Math. Soc. (2) 
25 (1982), no. 1, 25--34], {J. London Math. Soc.} (2) \textbf{46} (1992), no. 2, 202--204.

\bibitem{Bir1961}
Birch, B. J., Waring's problem in algebraic number fields, \emph{Proc. Cambridge Philos. Soc.} \textbf{57} 
(1961), 449--459.

\bibitem{Bok1994}
Boklan, K. D., The asymptotic formula in Waring's problem, \emph{Mathematika} \textbf{41} (1994), 
no. 2, 329--347.

\bibitem{BW2012}
Boklan, K. D. and Wooley, T. D., On Weyl sums for smaller exponents, \emph{Funct. Approx. Comment. 
Math.} \textbf{46} (2012), no. 1, 91--107.

\bibitem{BB2010}
Blomer, V. and Br\"udern, J., The number of integer points on Vinogradov's quadric, \emph{Monatsh. 
Math.} \textbf{160} (2010), no. 3, 243--256.

\bibitem{Bom1990}
Bombieri, E., On Vinogradov's mean value theorem and Weyl sums, \emph{Automorphic forms and 
analytic number theory (Montreal, PQ, 1989)}, pp. 7--24, Univ. Montr\'eal, Montreal, QC, 1990. 

\bibitem{BC2003}
Bourgain, J. and Chang, M.-C., On the size of $k$-fold sum and product sets of integers, 
\emph{J. Amer. Math. Soc.} \textbf{17} (2004), no. 2, 473--497.

\bibitem{Cha2003}
Chang, M.-C., The Erd\H os-Szemer\'edi problem on sum set and product set, \emph{Ann. of Math.} (2) 
\textbf{157} (2003), no. 3, 939--957.

\bibitem{CCG2014}
Chang, M.-C., Cilleruelo, J., Garaev, M. Z., Hern\'andez, J., Shparlinski, I. E. and Zumulac\'arregui, A., 
Points on curves in small boxes and applications, submitted, arXiv:1111.1543.

\bibitem{vdC1922}
van der Corput, J. G., Verscharfung der Absch\"atzungen beim Teilerproblem, \emph{Math. Ann.} 
\textbf{87} (1922), 39--65.

\bibitem{CH2010}
Croot, E. and Hart, D., $h$-fold sums from a set with few products, \emph{SIAM J. Discrete Math.} 
\textbf{24} (2010), no. 2, 505--519.

\bibitem{ES1983}
Erd\H os, P. and Szemer\'edi, E., On sums and products of integers, \emph{Studies in Pure Mathematics}, 
Birkh\"auser, Basel, 1983, pp. 213--218.

\bibitem{For1995}
Ford, K. B., New estimates for mean values of Weyl sums, \emph{Internat. Math. Res. Notices} (1995), 
no. 3, 155--171.

\bibitem{For2002}
Ford, K. B., Vinogradov's integral and bounds for the Riemann zeta function, \emph{Proc. London Math. 
Soc.} (3) \textbf{85} (2002), no. 3, 565--633.

\bibitem{FW2014}
Ford, K. B. and Wooley, T. D., On Vinogradov's mean value theorem: strongly diagonal behaviour via 
efficient congruencing, submitted, arXiv:1304.6917.

\bibitem{HL1922}
Hardy, G. H. and Littlewood, J. E., Some problems of `Partitio Numerorum': IV. The singular series in 
Waring's Problem and the value of the number $G(k)$, \emph{Math. Zeit.} \textbf{12} (1922), 161--188.

\bibitem{HB1988}
Heath-Brown, D. R., Weyl's inequality, Hua's inequality, and Waring's problem, \emph{J. London Math. 
Soc.} (2) \textbf{38} (1988), no. 2, 216--230.

\bibitem{Hil1909}
Hilbert, D., Beweis f\"ur die Darstellbarkeit der ganzen Zahlen durch eine feste Anzahl $n^{\text{ter}}$ 
Potenzen (Waringsches Problem), \emph{Math. Ann.} \textbf{67} (1909), no. 3, 281--300.

\bibitem{Hua1938}
Hua, L.-K., On Waring's problem, \emph{Quart. J. Math. Oxford} \textbf{9} (1938), 199--202.

\bibitem{Hua1938b}
Hua, L.-K., On Tarry's problem, \emph{Quart. J. Math. Oxford} \textbf{9} (1938), 315--320.

\bibitem{Hua1947}
Hua, L. K., \emph{The additive prime number theory}, Trav. Inst. Math. Stekloff, \textbf{22}, Acad. Sci. 
USSR, Moscow-Leningrad, 1947.

\bibitem{Hua1949}
Hua, L.-K., An improvement of Vinogradov's mean value theorem and several applications, 
\emph{Quart. J. Math. Oxford} \textbf{20} (1949), 48--61.

\bibitem{Hua1949b}
Hua, L.-K., Improvement of a result of Wright, \emph{J. London Math. Soc.} \textbf{24} (1949), 
157--159.

\bibitem{Kar1973}
Karatsuba, A. A., The mean value of the modulus of a trigonometric sum, \emph{Izv. Akad. Nauk SSSR 
Ser. Mat.} \textbf{37} (1973), 1203--1227.

\bibitem{Kor1962}
K\"orner, O., \"Uber Mittelwerte trigonometrischer Summen und ihre Anwendung in algebraischen 
Zahlk\"orpern, \emph{Math. Ann.} \textbf{147} (1962), 205--239.

\bibitem{Kor1958}
Korobov, N. M., Estimates of trigonometric sums and their applications, \emph{Uspehi Mat. Nauk} 
\textbf{13} (1958), no. 4 (82), 185--192.

\bibitem{KLZ2014}
Kuo, W., Liu, Y.-R. and Zhao, X., Multidimensional Vinogradov-type estimates in function fields, 
\emph{Canad. J. Math.}, in press.

\bibitem{Lin1943}
Linnik, Yu. V., On Weyl's sums, \emph{Mat. Sbornik (Rec. Math.) N. S.} \textbf{12} (1943), 28--39.

\bibitem{Mit1986}
Mit$^\prime$kin, D. A., Estimate for the number of summands in the Hilbert-Kamke problem, 
\emph{Mat. Sbornik (N.S.)} \textbf{129} (1986), no. 4, 549--577.

\bibitem{Mit1987}
Mit$^\prime$kin, D. A., Estimate for the number of summands in the Hilbert-Kamke problem, II, 
\emph{Mat. Sbornik (N.S.)} \textbf{132} (1987), no. 3, 345--351. 

\bibitem{Par2005}
Parsell, S. T., A generalization of Vinogradov's mean value theorem, \emph{Proc. London Math. Soc.} (3) 
\textbf{91} (2005), no. 1, 1--32.

\bibitem{Par2014}
Parsell, S. T., A note on Weyl's inequality for eighth powers, \emph{Rocky Mountain J. Math.}, in press.

\bibitem{PPW2013}
Parsell, S. T., Prendiville, S. M. and Wooley, T. D., Near-optimal mean value estimates for multidimensional 
Weyl sums, \emph{Geom. Funct. Anal.} \textbf{23} (2013), no. 6, 1962--2024. 

\bibitem{Pre2013}
Prendiville, S. M., Solution-free sets for sums of binary forms, \emph{Proc. London Math. Soc.} (3) 
\textbf{107} (2013), no. 2, 267--302. 

\bibitem{RS2000}
Robert, O. and Sargos, P., Un th\'eor\`eme de moyenne pour les sommes d'exponentielles. Application 
\`a l'in\'egalit\'e de Weyl, \emph{Publ. Inst. Math. (Beograd) (N.S.)} \textbf{67} (2000), 14--30.

\bibitem{Rog1986}
Rogovskaya, N. N., An asymptotic formula for the number of solutions of a system of equations, 
\emph{Diophantine Approximations}, Part II, Moskov. Gos. Univ., Moscow, 1986, pp. 78--84.

\bibitem{Ste1975}
Stechkin, S. B., On mean values of the modulus of a trigonometric sum, \emph{Trudy Mat. Inst. Steklov} 
\textbf{134} (1975), 283--309.

\bibitem{Tyr1987}
Tyrina, O. V., A new estimate for a trigonometric integral of I. M. Vinogradov, \emph{Izv. Akad. Nauk 
SSSR Ser. Mat.} \textbf{51} (1987), no. 2, 363--378.

\bibitem{Ust1998}
Ustinov, A. V., On the number of summands in the asymptotic formula for the number of solutions of the 
Waring equation, \emph{Mat. Zametki} \textbf{64} (1998), no. 2, 285--296.

\bibitem{Vau1986a}
Vaughan, R. C., On Waring's problem for cubes, \emph{J. Reine Angew. Math.} \textbf{365} (1986), 
122--170.

\bibitem{Vau1986b}
Vaughan, R. C., On Waring's problem for smaller exponents, II, \emph{Mathematika} \textbf{33} (1986), 
no. 1, 6--22.

\bibitem{Vau1997}
Vaughan, R. C., \emph{The Hardy-Littlewood method}, 2nd edition, Cambridge University Press, 
Cambridge, 1997.

\bibitem{VW1995}
Vaughan, R. C. and Wooley, T. D., On a certain nonary cubic form and related equations, 
\emph{Duke Math. J.} \textbf{80} (1995), no. 3, 669--735.

\bibitem{VW1997}
Vaughan, R. C. and Wooley, T. D., A special case of Vinogradov's mean value theorem, \emph{Acta Arith.} 
\textbf{79} (1997), no. 3, 193--204. 

\bibitem{Vin1935}
Vinogradov, I. M., New estimates for Weyl sums, \emph{Dokl. Akad. Nauk SSSR} \textbf{8} (1935), 
195--198.

\bibitem{Vin1947}
Vinogradov, I. M., The method of trigonometrical sums in the theory of numbers, \emph{Trav. Inst. Math. 
Stekloff} \textbf{23} (1947), 109pp.

\bibitem{Vin1958}
Vinogradov, I. M., A new estimate of the function $\zeta(1+it)$, \emph{Izv. Akad. Nauk SSSR. Ser. Mat.} 
\textbf{22} (1958), 161--164.

\bibitem{Wal1963}
Walfisz, A. Z., \emph{Weylsche Exponentialsummen in der neueren Zahlentheorie}, Deutscher Verlag der 
Wissenschaften, Berlin, 1963.

\bibitem{Wey1916}
Weyl, H., \"Uber die Gleichverteilung von Zahlen mod Eins, \emph{Math. Ann.} \textbf{77} (1916), 
313--352.

\bibitem{Woo1992a}
Wooley, T. D., Large improvements in Waring's problem, \emph{Ann. of Math.} (2) \textbf{135} (1992), 
no. 1, 131--164.

\bibitem{Woo1992b}
Wooley, T. D., On Vinogradov's mean value theorem, \emph{Mathematika} \textbf{39} (1992), no. 2, 
379--399.

\bibitem{Woo1993}
Wooley, T. D., On Vinogradov's mean value theorem, II, \emph{Michigan Math. J.} \textbf{40} (1993), 
no. 1, 175--180.

\bibitem{Woo1994}
Wooley, T. D., Quasi-diagonal behaviour in certain mean value theorems of additive number theory, 
\emph{J. Amer. Math. Soc.} \textbf{7} (1994), no. 1, 221--245. 

\bibitem{Woo1995}
Wooley, T. D., New estimates for Weyl sums, \emph{Quart. J. Math. Oxford} (2) \textbf{46} (1995), 
no. 1, 119--127.

\bibitem{Woo1996}
Wooley, T. D., Some remarks on Vinogradov's mean value theorem and Tarry's problem, \emph{Monatsh. 
Math.} \textbf{122} (1996), no. 3, 265--273.

\bibitem{Woo2002a}
Wooley, T. D., Diophantine methods for exponential sums, and exponential sums for Diophantine problems, 
{\it Proceedings of the International Congress of Mathematicians, August 20--28, 2002, Beijing}, Volume II, 
Higher Education Press, 2002, pp. 207--217.

\bibitem{Woo2012a}
Wooley, T. D., Vinogradov's mean value theorem via efficient congruencing, \emph{Ann. of Math.} (2) 
\textbf{175} (2012), no. 3, 1575--1627.

\bibitem{Woo2012b}
Wooley, T. D., The asymptotic formula in Waring's problem, \emph{Internat. Math. Res. Notices} (2012), 
no. 7, 1485--1504.

\bibitem{Woo2013}
Wooley, T. D., Vinogradov's mean value theorem via efficient congruencing, II, \emph{Duke Math. J.} 
\textbf{162} (2013), no. 4, 673--730.

\bibitem{Woo2014a}
Wooley, T. D., Multigrade efficient congruencing and Vinogradov's mean value theorem, submitted, 
arXiv:1310.8447.

\bibitem{Woo2014b}
Wooley, T. D., Approximating the main conjecture in Vinogradov's mean value theorem, submitted, 
arXiv:1401.2932.

\bibitem{Woo2014c}
Wooley, T. D., The cubic case of the main conjecture in Vinogradov's mean value theorem, submitted, arXiv:1401.3150.

\bibitem{Woo2014d}
Wooley, T. D., Mean value estimates for odd cubic Weyl sums, submitted, arXiv:1401.7152.

\bibitem{Wri1948}
Wright, E. M., The Prouhet-Lehmer problem, \emph{J. London Math. Soc.} \textbf{23} (1948), 279--285.

\end{thebibliography}
\end{document}